%% file: article.tex
\documentclass[a4paper]{amsart}

\usepackage{amssymb}
\usepackage{xypic}
\usepackage[pdfusetitle]{hyperref} 
\usepackage{tikz}
\usepackage{etoolbox}
\usepackage{subcaption}

\begin{document}

\newtheorem{theorem}[equation]{Theorem}
\newtheorem{conjecture}[equation]{Conjecture}
\newtheorem{proposition}[equation]{Proposition}
\newtheorem{lemma}[equation]{Lemma}
\newtheorem{corollary}[equation]{Corollary}
\theoremstyle{definition}
\newtheorem{definition}[equation]{Definition}
\newtheorem{remark}[equation]{Remark}
\newtheorem{example}[equation]{Example}
\theoremstyle{plain}
\numberwithin{equation}{section} 
\preto{\table}{\stepcounter{equation}}
\renewcommand\thetable{\thesection.\arabic{equation}}
\preto{\figure}{\stepcounter{equation}}
\renewcommand\thefigure{\thesection.\arabic{equation}}

\title[$p$-adic Artin $L$-functions]{Conjectures, consequences, and numerical\\ experiments for $p$-adic Artin $L$-functions}

\author{Rob de Jeu}
\address{Faculteit der B\`etawetenschappen\\Afdeling Wiskunde\\Vrije Universiteit Amsterdam\\De Boelelaan 1081a\\1081 HV Amsterdam\\The Netherlands}
\author{Xavier-Fran\c cois Roblot}
\address{Institut Camille Jordan \\ Universit\'e de Lyon, Universit\'e
  Lyon~1, CNRS -- UMR 5208 \\ 43 blvd du 11 Novembre 1918, 69622 Villeurbanne Cedex \\ France}

\begin{abstract}
We conjecture that the $ p $-adic $ L $-function of a non-trivial
irreducible even Artin character over a totally real field is non-zero
at all non-zero integers. This implies that a conjecture formulated
by Coates and Lichtenbaum at negative integers extends in a suitable way to all
positive integers.
We also state a conjecture that for certain characters the
Iwasawa series underlying the $ p $-adic~$ L $-series have no multiple roots except for those
corresponding to the zero at $ s = 0 $ of the $p$-adic $L$-function.

We provide some theoretical evidence for our first conjecture,
and prove both conjectures by means of computer calculations
for a large set of characters (and integers where appropriate) over the rationals and over real quadratic fields,
thus proving many instances of conjectures by Coates and Lichtenbaum and by Schneider.
The calculations and the theoretical evidence also prove that certain $ p $-adic regulators corresponding
to $ 1 $-dimensional characters for the rational numbers are units in many cases.
We also verify Gross' conjecture for the order of the zero of
the $ p $-adic $ L $-function at~$ s = 0 $ in many cases.
We gather substantial statistical data on the constant term of the underlying Iwasawa series,
and propose a model for its behaviour for certain characters.
\end{abstract}



\maketitle

\input defs.tex
\input intro.tex

\input Lp.tex

\input nozeroes.tex
\input consequences.tex

\input NP.tex

\input zeronumevidence.tex
\input lambdamodel.tex

\bibliographystyle{plain}
\bibliography{References}

\end{document}

%% file: defs.tex
\def\spaceifletter{\futurelet\comingchar\dospaceifletter}
\def\dospaceifletter{\relax\ifmmode\else
  \ifcat A\noexpand\comingchar{} \fi
  \ifcat 0\noexpand\comingchar
  \ifx 0\noexpand\comingchar{} \fi
  \ifx 1\noexpand\comingchar{} \fi\ifx 2\noexpand\comingchar{} \fi
  \ifx 3\noexpand\comingchar{} \fi\ifx 4\noexpand\comingchar{} \fi
  \ifx 5\noexpand\comingchar{} \fi\ifx 6\noexpand\comingchar{} \fi
  \ifx 7\noexpand\comingchar{} \fi\ifx 8\noexpand\comingchar{} \fi
  \ifx 9\noexpand\comingchar{} \fi\fi
  \ifcat $\noexpand\comingchar{} \fi
  \ifcat \noexpand\relax\noexpand\comingchar{} \fi
\fi}

\def\l{\lambda}
\def\o{\omega}

\def\C{\mathbb C}
\def\D{\mathcal D}
\def\et{\textup{\'et}}
\def\etale{\'etale\spaceifletter}

\def\Gal{\operatorname{Gal}}
\def\Ind #1 #2 {\textup{Ind}_{#1}^{#2}}
\def\Lfunction{$L$-function\spaceifletter}
\def\Lfunctions{$L$-functions\spaceifletter}
\def\Lp #1 #2 #3 #4 {L_{#1}(#2,#3,#4)}
\def\O{{\mathcal O}}
\def\Oe{\mathcal{O}_{E}}
\def\Q{\mathbb Q}
\def\Qbar{\overline{\mathbb Q}}
\def\Qp{{\mathbb Q_p}}
\def\Qpbar{{\overline {\mathbb Q}_p}}
\def\padic{$p$-adic\spaceifletter}
\def\phi{\varphi}
\def\rightiso{\buildrel{\sim}\over{\rightarrow}}
\def\tensor{\otimes}
\def\Z{\mathbb Z}
\def\Zp{{\mathbb Z}_p}

\def\Eul{\textup{Eul}}
\def\Eulpol#1{F_{#1}}
\def\dEul#1{\textup{Eul}_{#1}^{\langle\rangle}}
\def\Nm{\text{Nm}}
\def\Fr{\textup{Fr}}

\def\hchi{{\hat\chi}}

\def\ol{\overline}

\newcommand{\Art}{\mathrm{Art}}

%% file: intro.tex
\section{Introduction}

We fix an algebraic closure $ \Qbar $ of $ \Q $.  If $ k $ is  a subfield of $ \Qbar $
then we write $ G_k $ for the Galois group $ \Gal(\Qbar/k) $.
We view all number fields as subfields of $ \Qbar $,
hence we have $ G_k \subseteq G_\Q $.
We also fix a prime number $ p $ and an algebraic closure $ \Qpbar $
of~$ \Q_p $, with absolute value $ |\cdot|_p $ and valuation $v_p$
normalized by $ |p|_p = p^{-1} $ and $ v_p(p) = 1 $.

Let $ k $ be a number field,
$ \eta : G_k \to \Qpbar $ an Artin character.
We write $ \Qp(\eta) $ for the subfield of $ \Qpbar $ obtained
by adjoining all values of $ \eta $ to $ \Qp $, and $ \Z_{p,\eta} $
for the valuation ring of $ \Qp(\eta) $. 

Fix an integer $ m\le 0 $.
Let $\sigma: \Qpbar \to \C$ be an embedding.
If $ \eta $ is 1-dimensional, then it follows from \cite[VII Corollary 9.9]{Neu99}
that the value $L(k,\sigma \circ \eta, m)$ of the classical Artin
$ L $-function is in $\sigma(\Qpbar)\subseteq \C$ and that
\begin{equation*}
  L^\ast(k,\eta,m) = \sigma^{-1}(L(k,\sigma \circ \eta,m))
\end{equation*}
is independent of the choice of $\sigma$.  By the discussion
at the beginning of~\cite[\S~3]{dJNa} this also holds if~$ \eta $
is not 1-dimensional.
In particular, if $ \eta $ takes values in a finite extension
$ E $ of $ \Q $, then all those values are in $ E $.
Clearly,
if $ v $ is a finite place of $ k $, then
$ \Eul_v^\ast(k,\eta,m) = \sigma^{-1}(\Eul_v(k,\sigma \circ \eta,m)) $,
with $ \Eul_v(k,\sigma\circ\eta,s) $ the (reciprocal of the) Euler factor
for~$ v $, is independent of $ \sigma $ and lies in~$ E $ (see ~\eqref{recEulerfactor} 
for the expression of this Euler factor). 

Write $ q = p $ if $ p \ne 2 $, and $ q = 4 $ if $ p = 2 $.  
Let $ \o_p : G_k \to \Qpbar $ be the Teichm\"uller character
for $ p $, which is the composition
\begin{equation*}
 G_k \to G_\Q \to \Gal(\Q(\mu_{p^\infty})/\Q) \rightiso \Zp^* \to \mu_{\phi(q)}
\,,
\end{equation*}
where the last homomorphism maps
an element of $\Zp^*$ to the unique element of the torsion subgroup
$\mu_{\phi(q)} $ of $ \Zp^* $ in its coset modulo $1+q\Zp$.
Note that for $ k = \Q $, $ \o_p $ gives an isomorphism $ \Gal(\Q(\mu_q)/\Q) \simeq \Zp^*/(1+q\Zp) \simeq \mu_{\phi(q)} $.
We shall also use the projection  $ \Zp^* = \mu_{\phi(q)}  \times (1+q\Zp) \to 1+ q \Zp$,
mapping $x$ to $ \langle x \rangle$, so that $ x^{-1} \langle x \rangle $
is in~$ \mu_{\phi(q)} $.

Now assume that $ k $ is totally real. If $ E $ is any extension of $ \Q $,
and $ \eta : G_k \to E $ an Artin character of~$ k $, then $ \eta $ is called even if $ \eta(c) = \eta(1) $ for each $ c $
in the conjugacy class of complex conjugation in $ G_k $.  This is equivalent with the fixed field
$ k_\eta $ of the kernel of the underlying representation of $ \eta $ being totally real.
We call $ \eta $ odd if $ \eta(c) = -\eta(1) $ for those $ c $,
so that $ \eta $ is odd if and only if all these $ c $ act as multiplication by $ -1 $ in the underlying representation.

For $ \chi : G_k \to \Qpbar $ a 1-dimensional even Artin character
of the totally real number field $ k \subset \Qbar $, there is
a unique continuous function $ \Lp p k {\chi} s $ on $ \Zp $,
with values in $ \Qp(\chi) $, such that for all integers
$ m \le 0 $ we have the interpolation formula
\begin{equation} \label{interpol}
 \Lp p k {\chi} m = L^*(k , \chi\o_p^{m-1} , m) \prod_{v \in P}  \Eul_v^\ast(k,\chi\o_p^{m-1},m)
\,,
\end{equation}
where $ P $ consists of the places of $ k $ lying above~$ p $.
The right-hand side here is never zero if $ m < 0 $, and $ L^*(k , \chi\o_p^{-1} , 0) \ne 0 $
as well.

Such continuous functions in particular cases (see, e.g., 
\cite{Kub-Leo,Ser73}) were the starting point, but, in fact,
they have much nicer properties as they can also be described as follows
(see, e.g., \cite{Bar78, Cas-Nog79, Del-Rib, Rib79}, although we
shall mostly follow \cite{greenberg83} here).

Let $ r \ge 1 $ be the largest integer such that $ \mu_{p^r} $ is in $ k(\mu_q) $
and put
\begin{equation*}
 \D_k = \{ z \text{ in } \Qpbar \text{ with } |z|_p < p^r p^{-1/(p-1)} \} 
\,.
\end{equation*}
Fix a topological generator $\gamma$ of $ \Gal(k_\infty/k) $,
where $ k_\infty $ is the cyclotomic $ \Zp $-extension of $ k $.
(We shall briefly review in Remark~\ref{indeprem} how the various
notions in this introduction depend on this choice.)
Let $ \gamma' $ in $ \Gal(k_\infty(\mu_q)/k) $ be the unique
element that restricts to $ \gamma $ on $ k_\infty $ and to the identity
on $ k(\mu_q) $.
Since the composition
\begin{equation*}
 \Gal(k_\infty(\mu_q)/k) \to \Gal(\Q_\infty(\mu_q)/\Q) \rightiso \Z_p^* 
\end{equation*}
is injective, there is a unique $ u $ in $ 1+q\Zp $ such that $ \gamma'(\xi) = \xi^u $
for all roots of unity~$ \xi $ of $ p $-power order.

Let $H_\chi(T)= \chi(\gamma)(1+T) - 1$ 
if $k_\chi \subset k_\infty $ (i.e., $\chi$ is of type $W$) and $ H_\chi(T) = 1 $ otherwise.
There exists a unique power series $G_\chi(T)$, called the Iwasawa power series of $\chi$,
with coefficients in $\Z_{p,\chi}$, such that,
for all $s $ in $ \Zp $, with $s \ne 1$ if $\chi$ is trivial, we have
\begin{equation} \label{iwasawaseries}
 \Lp p k {\chi} {s} = \frac{G_\chi(u^{1-s} - 1)}{H_\chi(u^{1-s} - 1)}
\,.
\end{equation}
Note that $ | u-1 |_p = p^{-r} $, so $ u^t = \exp(t\log(u)) $
exists for all $ t $ in~$ \D_k $. Therefore $ \Lp p k {\chi} s $
is a meromorphic function on $ \D_k $, 
with values in $ \Qpbar $,
with at most a pole of
order~1 at $ s = 1 $ if $ \chi $ is trivial, and no poles otherwise.
We let $ \zeta_p(k,s) $ denote this function if $ \chi $ is the trivial character.

If $ \chi : G_k \to \Qpbar $ is an even irreducible Artin character
of degree greater than~1, we let $ H_\chi(T) = 1 $.
We then define $ H_\chi(T) $ for all characters $ \chi $
by demanding $ H_{\chi_1+\chi_2}(T) = H_{\chi_1}(T) H_{\chi_2}(T) $.
Using the Galois action on $ \Qpbar $ one sees that $ H_\chi(T) $ has coefficients
in $ \Z_{p,\chi} $.\footnote{In \cite[\S2]{greenberg83} 
a ring $ \Zp[\chi] $ is used for all $ \chi $,
but there it is only defined for $ \chi $ 1-dimensional, as
the ring  obtained by adjoing the values of $ \chi $ to $ \Zp $.
At the bottom of page 82 of loc.\ cit.\ it is clear that $ \Z_{p,\chi} $
is meant. Note that $ \Zp[\chi] = \Z_{p,\chi} $ if $ \chi $ is
1-dimensional, and that $ \Zp[\textup{values of } \chi] $ is contained
in $ \Z_{p,\chi} $ for all $ \chi $ but may be smaller even for Abelian $ \chi $.}

For $ \chi : G_k \to \Qpbar $ any even Artin character of $ G_k $, one
uses Brauer induction to find a unique $ G_\chi(T) $ in the fraction
field of $ \Z_{p,\chi}[[T]] $
such that~\eqref{interpol} holds for all integers $ m < 0 $ if
we define $ \Lp p k {\chi} s $ by~\eqref{iwasawaseries} using
this $ G_\chi(T) $ and the $ H_\chi(T) $ defined above.
(We shall briefly review this in Section~\ref{Lpsection}.)
It is known for $ p $ odd (and conjectured for $ p = 2 $) that
$ G_\chi $ is in $ p^{-l} \Z_{p,\chi}[[T]] $ for some integer $ l $,
in which case it converges on~$ \D_k $.
This means that $ \Lp p k {\chi} s $
is defined on $ \D_k $ for $ p $ odd, but for $ p=2 $ it might
have finitely many poles if $ \chi $ is not Abelian.

We can now state our main conjecture,
some special cases of which appeared in the literature before
(see Remark~\ref{earlier}).
The equivalence of the
various parts will be proved in Section~\ref{nozeroessection}.
The `missing' case $ m = 0 $ in this conjecture is the subject
of a conjecture by Gross, Conjecture~\ref{gross conj} below.
Note that the statement here always holds for $ m < 0 $ by~\eqref{interpol}
as the right-hand side there is never zero, but we included this
case in the formulation of the conjecture anyway for the sake of uniformity.

\begin{conjecture} \label{mainconjecture}
Fix a prime number $ p $ and an integer $ m \ne 0 $.
Then the following equivalent statements hold, where in the last
three parts the characters always take values in~$ \Qpbar $.
\begin{enumerate}
\item
For every totally real number field $ k $ with $ k/\Q $ Galois
the function $ \zeta_{k,p}(s) $ is non-zero at $ m $ if $ m \ne 1 $ and has a pole of order~1 at $ m = 1 $.

\item
For every totally real number field $ k $ the function $ \zeta_{k,p}(s) $
is non-zero at $ m $ if $ m \ne 1 $ and has a pole of order~1 at $ m = 1 $.

\item
For every totally real number field $ k $ and every 1-dimensional
even Artin character $ \chi \ne {\bf 1}_{G_k} $ of $ G_k $ we have $ \Lp p k {\chi} m \ne 0 $.

\item
For every totally real number field $ k $ and every irreducible
even Artin character $ \chi \ne {\bf 1}_k $ of $ G_k $ we have that $ \Lp p k {\chi} m $ is defined and non-zero.

\item
For every irreducible even Artin character
$ \chi \ne {\bf 1}_{\Q} $ of $ G_\Q $ we have that $ \Lp p {\Q} {\chi} m $ is defined and non-zero.
\end{enumerate}
\end{conjecture}

For $ p $ odd, one can state the conjecture also in terms of
the power series $ G_\chi(T) $ associated to $ k $ and $ \chi $ that we introduced above,
and for $ p=2 $ using fractions of power series. It is more natural
in the sense that $ m = 1 $ is no longer a special case; see Remark~\ref{Gremark}.
We prove this version in a number of cases (see Lemma~\ref{vlemma} and
Proposition~\ref{vprop}).

Our numerical calculations  prove Conjecture~\ref{mainconjecture}(3) for many tuples $ (p,k,\chi,m)$ with
$ [k : \Q ] = 1 $ or~2, and in fact, for most tuples $ (p,k,\chi) $
considered the statement holds for all $ m \ne 0 $; see Theorem~\ref{th:mainconj}.
Because Conjecture~\ref{mainconjecture}
has important consequences for
(generalisations of) conjectures made by Coates and Lichtenbaum and by Schneider
(see Section~\ref{consequences}), our calculations prove many
instances of these conjectures.
The conjecture
was itself inspired mostly by the wish
that certain $ p $-adic regulators in algebraic $ K $-theory
should be non-zero, and our calculations also prove this in many
cases (see the end of Section~\ref{consequences}). We can also prove it in certain
cases without relying on computer calculations (see Corollary~\ref{regcor}).

\smallskip

In order to state our second conjecture we need more notation.

\smallskip

If $ \O $ is the valuation ring in a finite extension of $ \Qp $,
then according to~\cite[Theorem~7.3]{wash} any non-zero $ S(T) $
in $ \O[[T]] $ can be written uniquely as $ c D(T) U(T) $,
with $ c \ne 0 $ in $ \O $, $ U(T) $ a unit of $ \O[[T]] $ with $ U(0)=1 $, and $ D(T) $
a distinguished polynomial, i.e., a monic polynomial in $ \O(T) $
such that all non-leading coefficients are in the valuation ideal
of~$ \O $. If we extend the field then this decomposition remains
the same.
For $ S(T) = c D(T) U(T) $ as above we define $ \mu(S) = v_p(c) $
and we let $ \l(S) $ be the degree of $ D(T) $.
We extend those definitions to the field of fractions of $ \O[[T]] $
in the obvious way.

\newpage

We can now state our second conjecture.

\begin{conjecture} \label{Grootsconjecture}
Let $\chi : G_k \to \Qpbar $ be a 1-dimensional even Artin character
of a totally real number field~$ k $.
Assume that $\Ind {G_k} {G_\Q} (\chi)$ is irreducible.
Then any root of the distinguished polynomial $ D_\chi(T) $ of $ G_\chi(T) $ distinct from $u-1$ is simple.
\end{conjecture}

Note that in~\eqref{iwasawaseries} we have $ |u^{1-s}-1|_p < p^{-1/(p-1)} $,
so this conjecture is also about roots of $ D_\chi(T) $ that
are not detected by the zeroes of $ \Lp p k {\chi} s $ on~$ \D_k $.
The cases in which we prove this conjecture numerically are described
in Theorem~\ref{th:simple}.

The potential root $ u-1 $ of $ G_\chi(T) $ corresponds to the
potential root~0 of $ \Lp p k {\chi} s $. For
this root, Gross formulated the following as part of Conjecture~2.12 in~\cite{Gross81}.

\begin{conjecture} \label{gross conj}
Let $\chi : G_k \to \Qpbar $ be a 1-dimensional even Artin character
of a totally real number field~$ k $,
$ P $ the set of places of  $ k $ lying above~$ p $,
and $\sigma : \Qpbar \to \mathbb{C}$ any embedding.  Then
the $p$-adic $L$-function $\Lp p k {\chi} s $ 
and the truncated complex $L$-function 
$  L(k , \sigma \circ\chi\o_p^{-1} , s) \prod_{v \in P}  \Eul_v(k,\sigma \circ \chi\o_p^{-1},s) $
have the same order of vanishing at $ s=0 $.
\end{conjecture}

\begin{remark}
One can also consider the case of truncated \padic \Lfunctions
for a finite set~$ S $ of places of $ k $ that contains~$ P $,
and compare the order of vanishing at $ s = 0 $ to the corresponding
truncated complex \Lfunction, but, as noted in loc.~cit., the correctness of the
conjecture is independent of~$ S $.
\end{remark}

Because of the construction of $ \Lp p k {\chi} s $
for arbitrary $ \chi $, and the behaviour of the classical \Lfunctions
for Brauer induction, the same conjecture would then hold for
all even Artin characters $ \chi $ of~$ k $.
Theoretical evidence for it  will be discussed in Remark~\ref{zeroremark},
but we also prove this conjecture in many cases by means of our
calculations; see Theorem~\ref{th:proofgoss=0} and Remark~\ref{grossremark}.

We now formulate the last conjecture that our calculations prove in many cases. 
For an even Artin character $ \chi $ of $ G_k $, we define
$ \lambda(\chi) = \lambda(G_\chi/H_\chi) $ and
$ \mu(\chi) = \mu(G_\chi/H_\chi) $.
(Note that $ H_\chi $ is often ignored when defining $ \lambda(\chi) $
and $ \mu(\chi) $, but we follow \cite{Sinnott}.)
Clearly, $ \mu(H_\chi) = 0 $,
so $\mu(\chi) \geq 0$. Deligne and Ribet \cite{Del-Rib} proved
(see \cite[(4.8) and (4.9)]{Rib79})
that $\mu(\chi) \geq [k:\Q]$ if $p=2$.
The so-called ``$\mu = 0$'' conjecture, formulated in~\cite{Iwa73b}, states this
bound should be sharp, and that $ \mu(\chi) $ should be~0 if $ p $ is odd.

\begin{conjecture}[``$\mu = 0$''] \label{muconjecture}
Let $\chi : G_k \to \Qpbar $ be a 1-dimensional even Artin character
of a totally real number field~$ k $.
Then $\mu(\chi) = 0$ if $p$ is odd, and $\mu(\chi) = [k:\Q]$ if $p = 2$. 
\end{conjecture}

This conjecture was proved for $k = \Q$ by Ferrero and Washington \cite{Fer-Wash}.
Theorem~\ref{th:proofmu=0} describes for which real quadratic
number fields $ k $ and $ \chi $ our calculations prove it for
more than 6000 conjugacy classes over $\Qp$ of \padic
characters not covered by the result of Ferrero and Washington.

There do not seem to be many conjectures or results on the $ \lambda $-invariant
(but see \cite[Conjecture~1.3]{EJV11} for a very special case).
Based on the (extensive) data provided by our calculations we formulate
the following conjecture.

\begin{conjecture} \label{Qconjecture}
Let $p$ be an odd prime and let $d \geq 1$. Define $\mathfrak{X}_d$ to be the set of even $1$-dimensional Artin characters $\chi : G_\Q \to \Qpbar $ of the form~$\chi = \o_p^i\psi$ with $\psi$ of conductor and order prime to $p$, both $\psi$ and $i$ even, and such that~$[\Q_p(\chi):\Q_p] = d$. For $N \geq 1$, let $\mathfrak{X}_d(N)$ be the subset of those characters in $\mathfrak{X}_d$ whose conductor is at
most~$N$. Then we have 
\begin{equation*}
\lim_{N \to +\infty} \frac{\#\{\chi \in \mathfrak{X}_d(N) : \lambda(\chi) > 0\}}{\#\mathfrak{X}_d(N)} = p^{-d}. 
\end{equation*}
\end{conjecture}

During the calculations we discovered this conjecture does not apply
if we allow~$ p $ to divide the order of the character. In fact, in
that case the $ \lambda $-invariant can be quite large (see
Examples~\ref{firstexample} through~\ref{lastexample}).  This phenomenon
can be explained by a result of Sinott~\cite[Theorem~2.1]{Sinnott}; see Corollary~\ref{lambdappowercase}.
The conjecture does also not apply if we allow both
$ i $ and $ \psi $ to be odd; see the beginning of Section~\ref{lambdamodel}.

\medskip

The structure of this paper is as follows.

In Section~\ref{Lpsection} we
review and slightly extend the description of $ \Lp p k {\chi} s $ if
$ \chi $ is an even $ \Qpbar $-valued Artin character of higher
dimension.  We also review truncated \padic \Lfunctions, which play a
role in some of the conjectures or their consequences.

Section~\ref{nozeroessection} proves the equivalence of the various
parts of Conjecture~\ref{mainconjecture}, and gives a formulation
using the power series~$ G_\chi(T) $. It also discusses some earlier cases
of this conjecture and some theoretical evidence for it.
We conclude this section with a review of
what is known about Conjecture~\ref{gross conj},
and a discussion of our motivation for making Conjecture~\ref{Grootsconjecture}.

In Section~\ref{consequences} we discuss some important consequences for a
generalisation of a conjecture by Coates and Lichtenbaum
(\cite[Conjecture~1]{Co-Li}, but see also \cite[\S1]{dJNa}) and
for a conjecture by Schneider \cite{schn79}.
We also discuss the relation with a conjecture
made in~\cite{BBDJR}. In fact, this last conjecture inspired the
current paper as it is about the non-vanishing of \padic regulators in $ K $-theory,
in analogy with the Leopoldt conjecture, and our calculations
that are described in Section~\ref{zeronumevidence} prove this
non-vanishing in many cases.

In Section~\ref{NPsection} we review and slightly refine the theory of
Newton polygons, in order to help us rule out by computations
the existence of multiple factors in the corresponding distinguished polynomials in Section~\ref{zeronumevidence},
thus proving Conjecture~\ref{Grootsconjecture} in many cases.
In that section we also prove, by means of calculations, Conjectures~\ref{mainconjecture},
\ref{gross conj} and~\ref{muconjecture} for many characters.

Finally, in Section~\ref{lambdamodel}, we discuss the behaviour
of the $ \l $-invariants of the Iwasawa power series, and our
(substantial) numerical data leads us to make and corroborate Conjecture~\ref{Qconjecture}.

%% file: Lp.tex
\section{Review of $ p $-adic $ L $-functions} \label{Lpsection}

In this section, we collect from the literature some results 
that we shall need later on, occasionally clarifying or slightly
extending them. 

We begin with the definition of $ \Lp p k {\chi} s $
for $ \chi : G_k \to \Qpbar $ any even Artin character
of a totally real number field $ k $, and the construction of
 $ G_\chi $ and $ H_\chi $ satisfying~\eqref{iwasawaseries}.
This is based on \cite[\S2]{greenberg83}, with some minor extension
and correction, and some input from~\cite[\S3]{dJNa}.

Recall that Brauer's theorem (see, e.g., \cite[XVIII, Theorem~10.13]{lang93})
states that if $ G $ is a finite group and $ \chi $ a character of $ G $, then 
\begin{equation*}
 \chi = \sum_i a_i \Ind H_i G (\chi_i)
\end{equation*}
for some non-zero integers $ a_i $ and 1-dimensional characters $ \chi_i $ of subgroups $ H_i $
with each $ H_i $ the product of a cyclic group and a  group
of order a power of a prime~$ p_i $.

For future reference we observe that we may assume that $ \chi_i \ne {\bf 1}_{H_i} $ if $ H_i \ne \{e\} $,
in which case the multiplicity of $ {\bf 1}_G $ in $ \chi $ equals the coefficient of
$ \Ind \{e\} G ({\bf 1}_{\{e\}}) $ since by Frobenius reciprocity 
$ \Ind H_i G (\chi_i) $ does not contain $ {\bf 1}_G $ when $ \chi_i \ne {\bf 1}_{H_i} $.
Namely, every subgroup of an $ H_i $ is of the same type as we
may assume the cyclic group has order prime to~$ p_i $.
Also, $ [H_i,H_i] \ne H_i $ if $ H_i $ is non-trivial.
For a finite group $ H $ we have
$ \Ind [H,H] H ({\bf 1}_{[H,H]}) = {\bf 1}_H + \sum_j \chi_j  $ where the sum runs through the
non-trivial 1-dimensional characters of $ H $.  Using this and the transitivity of induction,
we can successively eliminate all ${\bf 1}_{H_i} $ 
unless $ H_i $ is trivial.

Applying the above to $ G = \Gal(k_\chi/k) $ we see that there exist
fields $ k_1,\dots,k_t $ with $ k \subseteq k_i \subseteq k_\chi $, 
1-dimensional Artin characters $\chi_i$ on $ G_{k_i} $, and non-zero integers $a_1, \dots, a_t$, with
\begin{equation} \label{induction}
  \chi = \sum_{i=1}^t a_i \Ind {G_{k_i}} {G_k} (\chi_i)
\,.
\end{equation}
If $ \chi $ does not contain the trivial character then 
we may assume that all $ \chi_i $ are non-trivial.
We then define the \padic \Lfunction of $\chi$ by
\begin{equation}\label{artinLp}
  \Lp p k {\chi} s = \prod_{i=1}^t \Lp p k_i {\chi_i} s ^{a_i}
\,.
\end{equation}
From~\eqref{iwasawaseries} one sees that it is a meromorphic function on $\D_k$
(cf., e.g., \cite[p.82]{greenberg83}, \cite[Section~6]{BBDJR}, \cite[\S3]{dJNa}).
If $ m $ is a negative integer satisfying $m \equiv 1 $ modulo $ {\varphi(q)} $
then the value $ \Lp p k {\chi} m $ is defined and equals
$ \Eul^\ast_p(k,\chi,m) L^*(k , {\chi} , m ) $
by~\eqref{interpol} and well-known properties of Artin
\Lfunctions~(see \cite[Prop.~VII.10.4(iv)]{Neu99}),
showing that the function is independent of how we express $\chi$ as a sum of induced 1-dimensional characters.
We say that $ \Lp p k {\chi} s $ is defined at $ s=s_0 $ if it does not have a pole at $ s=s_0 $, so
that we can talk about the value~$ \Lp p k {\chi} s_0 $.

If the index of $ \Gal(k_{i,\infty}/k_i) $ in $ \Gal(k_\infty/k) $ (both inside $ \Gal(\Q_\infty/\Q) $)
is $ p^{a_i} $, and we use the generator $ \gamma_i $ in the first
that corresponds to $ \gamma^{p^a_i} $, then $ u_i = u^{p^{a_i}} $.
If we define $ G_{\chi_i}^*(T) $ by substituting $ T_i = (T+1)^{p^{a_i}} - 1 $
in $ G_{\chi_i}(T_i) $, and similarly for~$ G_{\chi_i}^*(T) $,
then we have
\begin{equation} \label{twoTs}
    \Lp p k {\Ind G_{k_i} G_k (\chi_i)} s
  = \Lp p k_i {\chi_i} s 
  = \frac{G_{\chi_i}(u_i^{1-s}-1)}{H_{\chi_i}(u_i^{1-s}-1)}
  = \frac{G_{\chi_i}^*(u^{1-s}-1)}{H_{\chi_i}^*(u^{1-s}-1)} 
\end{equation}
for $ s $ in $ \D_k $.
Using this in~\eqref{artinLp} gives an expression as in~\eqref{iwasawaseries}
for some $ G_\chi(T) / H_\chi(T) $ in the fraction field of $ \Z_{p,\chi}[[T]] $,
as one can see using the action of $ \Gal(\Qpbar/\Qp) $ and the
fact that~\eqref{interpol} provides infinitely many $ u^{1-m}-1 $ in $ q\Zp $
where $ G_\chi(T) $ takes values in $ \Qp(\chi) $.
But, as proved by Greenberg, one can do much
better, based on the main conjecture of Iwasawa theory.

In the introduction we defined $ H_\chi(T) = 1 $ for $ \chi $ irreducible of degree
greater than~1, and $ H_\chi(T) $ for all characters $ \chi $
by demanding $ H_{\chi_1+\chi_2}(T) = H_{\chi_1}(T) H_{\chi_2}(T) $.
Then $ H_\chi(T) $ has coefficients in $ \Z_{p,\chi} $.

We want to compare this to Brauer induction.
Let $ k'/k $ be a finite extension with $ k' $ totally real,
$ \chi' $ an irreducible Artin character of $ k' $, and $ \chi_*' = \Ind G_{k'} G_k (\chi') $
the induced character on $ G_k $.
If we write down $ H_{\chi'}(T') $ for the generator $ \gamma' $
of $ \Gal(k_\infty'/k') $ corresponding to $ \gamma^{p^a} $ (again viewing everything in
$ \Gal(\Q_\infty/\Q) $), and let $ H_{\chi'}^*(T) = H_{\chi'}((T+1)^{p^a}  - 1 ) $,
then one checks that $ H_{\chi_*'}(T) $ and $ H_{\chi'}^*(T) $
are identically~1 unless $ \chi' $ is of type~$ W $, in which
case $ H_{\chi_*'}(T) = (-1)^{p^a+1} H_{\chi'}^*(T) $.
With $ G_{\chi'}^*(T) = G_{\chi'}((T+1)^{p^a}-1) $ we have
$ \Lp p k' \chi' s = G_{\chi'}^*(u^{1-s})/H_{\chi'}^*(u^{1-s}) $
for $ s $ in $ \D_k $.
Therefore from~\eqref{induction} one obtains~\eqref{iwasawaseries} for
all but finitely many $ s $ in $ \D_k $, if we let
the polynomial $ H_\chi(T) $ be as defined above, and
\begin{equation} \label{Ggeneral}
 G_\chi(T) = \pm \prod_{i=1}^t G_{\chi_i}^*(T)^{a_i}
\,,
\end{equation}
with the sign matching that in $ H_\chi(T) = \pm \prod_{i=1}^t H_{\chi_i}^*(T)^{a_i} $.
Again, $ G_\chi(T) $ is in the fraction field of $  \Z_{p,\chi}[[T]] $.

\begin{remark}
Greenberg claims that, given~\eqref{induction},
putting $ H_\chi(T) = \prod_{i=1}^t H_{\chi_i}^*(T)^{a_i} $ is
well-defined (see \cite[p.82]{greenberg83}).
But it follows from the above that this may fail if $ p = 2 $.
If $ k'/k $ is the extension of degree~2 in the cyclotomic
$ \Z_2 $-extension of $ k $, then using $ \chi = {\bf 1}_{G_k} $ or $ \chi = \Ind G_{k'} G_k ({\bf 1}_{G_{k'}}) - \chi' $
with $ \chi' $ the non-trivial character of $ G_k $ that is trivial
on $ G_{k'} $ leads to two different $ H_\chi(T) $ that differ by
a sign.

Of course, if one were to normalize all $ H_\chi(T) $ to be
monic, and absorb the resulting root of unity into $ G_\chi(T) $, this problem would disappear.
Here we opted to define $ H_\chi(T) $ uniquely, and with that
also $ G_\chi(T) $.
\end{remark}

One has the following conjecture about $ G_\chi(T) $ (see~\cite[\S2]{greenberg83}).

\begin{conjecture}{$(p $-adic Artin conjecture$)$} \label{artinconjecture}
Let $ k $ be a totally real number field, and~$ \chi : G_k \to \Qpbar $
an even \padic Artin character. Then $ p^l G_\chi(T) $ is in~$ \Z_{p,\chi}[[T]] $ for some~$ l $.
\end{conjecture}

Greenberg proved in \cite[Proposition~5]{greenberg83} that this
conjecture is implied by the main conjecture of Iwasawa theory.
Because the latter was proved for $ p $ odd by
Wiles in~\cite{Wiles90}, Conjecture~\ref{artinconjecture} holds for $ p \ne 2 $.
It then follows from~\eqref{iwasawaseries} that, for a character~$ \chi $,
and~$ p $ odd, $ \Lp p k {\chi} s $ is analytic on
$\D_k$ if $ \chi $ does not contain the trivial character, and has at most a pole at $ s=1 $ otherwise.
(This statement also followed from \cite{Wiles90} in a much more roundabout way in~\cite{dJNa}.)

\begin{remark}
Conjecture~\ref{artinconjecture} is equivalent with the quotient of distinguished polynomials
$  \prod_{i=1}^t D_{\chi_i}^*(T)^{a_i} $ 
corresponding to~\eqref{Ggeneral} being a polynomial. If this
is the case, then the statement of the conjecture holds precisely
when $ l + \mu(\chi) \ge 0 $.

Note that Conjecture~\ref{muconjecture} would imply that, for
any even Artin character $ \chi $ of $ G_k $, we have that $ \mu(\chi)  = 0 $
if $ p $ is odd, and $ \mu(\chi) = \chi(1) \cdot [k : \Q] $ if
$ p = 2 $.
(The latter statement one finds by applying~\eqref{induction}
to the element 1 of $ G_k $, multiplying the result by $ [k:\Q] $, and using that,
for an Artin character $ \psi $ on $ G_{k'} $, $ k \subseteq k' $,
and  $ \psi_* = \Ind G_{k'} G_k (\psi) $, one has
$ \psi(1) \cdot [k':\Q] = \psi_*(1) \cdot [k : \Q] $;
cf.\ the argument at the end of \cite[\S4]{greenberg83}.)
\end{remark}

\begin{remark} \label{indeprem}
As promised in the introduction we briefly explain that
the various statements there are independent of
the choice of $ \gamma $.
Namely, using $ \tilde\gamma = \gamma^a $ instead of $ \gamma $, with $ a $ in $ \Zp^* $,
replaces $ u $ with $ \tilde u = u^a  $.
Both $ G_\chi(T)/H_\chi(T) $ and $ \tilde G_\chi(\tilde T)/\tilde H_\chi(\tilde T) $
satisfy~\eqref{iwasawaseries} if we obtain the former from the
latter by replacing $ \tilde T $ with $ g(T) = (1+T)^a -1 = \sum_{n=1}^\infty \binom an T^n $
in $ T \Zp[[T]]^* $. 
Clearly, this substitution (as well as its inverse) maps units to units, and a distinguished
polynomial to a distinguished polynomial (of the same degree) times a unit, hence
they preserve the $ \l $-invariant (and trivially the $ \mu $-invariant).

An easy calculation (only needed when $ \chi $ is of type~$ W $) shows that
the distinguished polynomial of $ \tilde H_\chi(g(T)) $ equals that of $ H_\chi(T) $, 
so the corresponding statement holds for $ \tilde G_\chi(g(T)) $ and $ G_\chi(T) $.
In particular, they also have the same $ \l $-invariants, and matching
roots (with $ \tilde u -1 $ corresponding to $ u-1 $) together
with their multiplicities.
Expanding the right-hand side of~\eqref{iwasawaseries}, with
$ t = 1-s $, as element of $ \Qpbar((t)) $, one sees that the
extension to $ \D_k $ of $ \Lp p k {\chi} {s} $ is also independent
of the choice of $ \gamma $ as its value at infinitely many values
of $ t $ are fixed by~\eqref{interpol}.
\end{remark}

We conclude this section with a discussion of truncated \padic
\Lfunctions, which will be needed in Section~\ref{consequences}. For the sake of convenience, we shall follow~\cite[\S3]{dJNa}.

Let $ k $ be a number field, and $ \eta : G_k \to \Qpbar $ an Artin character.
Fix a finite place~$ v $ of $ k $. 
Let $ V $ be an Artin representation of $ G_k $, defined over a finite extension of~$ \Q(\eta) $,
with character $ \eta $.
If $ D_w $ is the decomposition group in $G_k$ of a prime $w$ lying above $v$,
with inertia subgroup $I_w$ and a Frobenius $ \Fr_w $, then we
let $ \Eulpol{v}( t, \eta) $ be the determinant of $ 1-\Fr_{w} t $
acting on $ V^{I_w} $ (note this is omitted just before \cite[(3.4)]{dJNa}, 
but mentioned just before \cite[(3.10)]{dJNa}).
Then $ \Eulpol{v} $ has coefficients in $ \Z_{p,\eta} $.
We put
\begin{equation} \label{recEulerfactor}
 \Eul_v^\ast(k,\eta,m) = \Eulpol{v}(\Nm(v)^{-m},\eta) 
\,.
\end{equation}
We can also use $ V $ to realize the character
$ \eta\o_p^{-1} $, and if $ v $ does not lie above
$ p $, then $ V^{I_w} $ is the same for both representations
as $ \o_p $ is unramified at~$ v $.
Clearly  $ \Eulpol{v}( \o_p(\Fr_w)^{m} t, \eta \o_p^{-m} ) = \Eulpol{v}( t, \eta) $
for any integer~$ m $.

Now assume that $ v $ does not lie above~$ p $.
If we write $ \Nm(v) $ for the norm of~$ v $, then
$ \Fr_w $ acts on $ p $-power roots of unity by raising
them to the power $ \Nm(v) $.
In $ \Gal(\Q(\mu_{p^\infty})/\Q) \simeq \Zp^* $, we can therefore
write  $ \Nm(v) = \langle \Nm(v) \rangle \cdot \o_p(\Fr_w) $
with $ \langle \Nm(v) \rangle $ in $ 1+q\Zp $.
For $ s $ in $ \Qpbar $ such that $ | s |_p \cdot |\langle \Nm(v) \rangle  -1 |_p <  p^{-1/(p-1)} $
we define
\begin{equation*}
 \dEul{v}(k,\eta,s) = \Eulpol{v}(\langle \Nm(v)\rangle^{-s},\eta) 
\,.
\end{equation*}
Then we have, for any integer~$ m $, that
\begin{equation*}
 \dEul{v}(k,\eta \o_p^{-m} ,m) =  \Eul_v^\ast(k,\eta,m)
\end{equation*}
because $ \Eulpol{v}(\langle \Nm(v) \rangle^{-m} ,\eta\o_p^{-m}) = \Eulpol{v}( \Nm(v)^{-m} ,\eta) $.
Note that, if $ \mu_{p^r} \subseteq k(\mu_q) $, then $ \langle \Nm(v) \rangle $ is
in $ 1 + p^r \Zp $ if $ p $ is odd or $  p = 2 $ and $ k $ is
totally real, because $ \Gal(k_{\mu_{p^\infty}} / k ) $, when identified
with a subgroup of $ \Gal(k_\infty/k) \times \Gal(k(\mu_q)/k) \subseteq (1+q\Zp) \times \mu_{\phi(q)} $,
is the product of subgroups of each of the factors.
In particular, in these cases this includes all $ s $ in $ \D_k $.

Now assume that $ k $ is totally real, and let $ \chi $ be an
even Artin character of $ G_k $. Let $ S $ be a finite set of
finite places of~$ k $, containing the set $ P $ of places above~$ p $.
For $ s $ in the domain of $ \Lp p k {\chi} s $
we define
\begin{equation} \label{pLS}
 \Lp {p,S} k {\chi} s = \Lp p k {\chi} s \prod_{v \in S \setminus P} \dEul{v}(k,\chi\o_p^{-1},s	)
\,.
\end{equation}
Note that $  \Lp p k {\chi} s =  \Lp {p,P} k {\chi} s $.
For an integer $ m<0 $ (and possibly for $ m=0 $),
$ \Lp {p,S} k {\chi} s $ is defined at $ m $, and in $ \Qp(\chi) $ we have
\begin{equation} \label{interpolS}
 \Lp {p,S} k {\chi} m =  L_S^\ast(k,\chi\o_p^{m-1},m) 
\end{equation}
where we let
$ L_S^\ast(k,\chi\o_p^{m-1},m) = L^\ast(k,\chi\o_p^{m-1},m) \prod_{v\in S} \Eul_v^\ast(k,\chi\o_p^{m-1},m) $.
Note that this also equals 
$ \sigma^{-1}( L_S(k,\sigma \circ \chi\o_p^{m-1},m) ) $ for any
embedding $ \sigma : \Qpbar \to \C $, because
$ \Eul_v^\ast(k,\eta,m) = \sigma^{-1}(\Eul_v(k,\sigma \circ \eta,m)) $
for any Artin character $ \eta : G_k \to \Qpbar $ and any integer~$ m $.
We also note that every $ \dEul{v}(k,\eta,s) $ in~\eqref{pLS}
is defined for $ s $ in $ \D_k $ and non-zero when $ s \ne 0 $ by \cite[Remark~3.11]{dJNa}.

%% file: nozeroes.tex
\section{Zeroes of $ p $-adic $ L $-functions} \label{nozeroessection}

In this section we discuss our main conjecture, Conjecture~\ref{mainconjecture}.

We start with proving that the parts of Conjecture~\ref{mainconjecture} are equivalent. We recall that
the Aramata-Brauer theorem 
(see \cite[Theorem~15.31]{CuRe1981} or \cite[XVIII, Theorem~8.4]{lang93})
 states the following.
If $ G $ is a finite group of order $ n>1 $ then for the character~$ \rho $ of the left-regular
representation we have $ n( \rho - 1 ) = \sum_i a_i \Ind H_i G (\chi_i) $
for some non-trivial 1-dimensional characters $ \chi_i $ of cyclic subgroups $ H_i $ of $ G $ and
positive integers $ a_i $.
In particular, if $ F/k $ is a Galois extension of totally real number fields of degree $ n \ge 1 $ then
either $ F=k $ or
\begin{equation}\label{ABcor}
\zeta_{F,p}(s)^n = \zeta_{k,p}(s)^n \prod_i \Lp p k_i \chi_i s ^{a_i}
\end{equation}
for some positive integers $ a_i $ and non-trivial 1-dimensional 
characters $ \chi_i $ of cyclic subgroups $ \Gal(F/k_i) $ of $ \Gal(F/k) $.

Returning to Conjecture~\ref{mainconjecture},
in order to see  that (1) implies~(2), let $ k $ be totally real and $ F $ its Galois closure over $ \Q $.  
If $ F=k $ then we are in case~(1).
If $ F \ne k $ then we notice that
in~\eqref{ABcor} the $ \Lp p k_i \chi_i s $ have no poles at all.
Because $ \zeta_{k,p}(s) $ can only have a pole (of order at
most~1) at $ s=1 $, and $ \zeta_{F,p}(s) $ by case~(1) has one
pole (of order~1) at $ s=1 $, $ \zeta_{k,p}(s) $ has one pole
(of order~1) at $ s=1 $.  At all integers $ m \ne 1 $ all functions
in~\eqref{ABcor} are defined, so $ \zeta_{F,p}(m) \ne 0 $ for
$ m \ne 0 $ implies $ \zeta_{k,p}(m) \ne 0 $ as well.

That~(2) implies~(3) is easily seen from $ \zeta_{F,p}(s) = \zeta_{k,p}(s) \prod_j \Lp p k \chi_j s $
where $ F $ is the fixed field of $ \ker(\chi) $ and the $ \chi_j $ run through the non-trivial 
1-dimensional characters of $ \Gal(F/k) $.

Deducing~(4) from~(3) is immediate since we may assume the $ \chi_j $ in~\eqref{artinLp} are non-trivial,
so that $ \Lp p k_j {\chi_j} s $ has no pole at $ s=m $, nor, by assumption, a zero.

That (4) implies (5) is clear, and that~(5) implies~(1) is clear if we write, for $ k/\Q $ Galois and
$ k $ totally real, $ \zeta_{k,p}(s) = \zeta_{\Q,p}(s) \prod_j \Lp p {\Q} \chi_j s ^{\chi_j(e)} $
where $ \chi_j $ runs through the non-trivial irreducible characters of $ G_\Q $ that are trivial
on $ G_k $: those $ \chi_j $ are necessarily even, and $ \zeta_{\Q,p}(s) $ has a pole of order 1 at
$ s=1 $ and no zeroes by~\cite[Lemma 7.12]{wash}.

This finished the proof of the equivalence of the statements
in Conjecture~\ref{mainconjecture}.

\begin{remark} \label{Gremark}
In the cases (1), (2) and~(3) of  Conjecture~\ref{mainconjecture},
because of~\eqref{iwasawaseries}, one can,
for a fixed $ k $ and $ \chi $, state the conjectured
behaviour of $ \Lp p k {\chi} s $ at $ m \ne 0 $ as
\begin{equation} \label{Gformulation}
  G_\chi(u^{1-m}-1) \ne 0
\,.
\end{equation}
Here $ \chi = {\bf 1}_{G_k} $ in (1) and (2), and $ \chi \ne {\bf 1}_{G_k} $ is
1-dimensional in~(3).

In parts (4) and (5), the same holds  if $ p \ne 2 $ because of
the discussion in Section~2, in particular, the part right after Conjecture~\ref{artinconjecture},
as $ G_\chi(u^{1-m}-1) $ is always defined.
For $ p = 2 $, we see from~\eqref{Ggeneral} that~\eqref{Gformulation} is defined and non-zero for $ \chi $,
is implied by $ G_{\chi_i}(u_i^{1-m}-1) = G_{\chi_i}^*(u^{1-m}-1) $
being non-zero for each of the $ \chi_i $ in~\eqref{induction}.

In particular, the conjecture for a fixed $ m \ne 0 $ is equivalent
with $ G_\chi(u^{1-m}-1) $ being non-zero for all even Artin characters~$ \chi $
of all totally real number fields~$ k $.
\end{remark}

We now discuss some known cases of Conjecture~\ref{mainconjecture}.
(numerical evidence for it when $ m \ge 1 $ will be discussed in Section~\ref{zeronumevidence}).
Note that it holds for all $ m < 0 $ as the right-hand side of~\eqref{interpol} is non-zero.

\medskip

We first consider the conjecture for $ \Lp p {\Q} \o_p^a s $
where $ a $ is an even integer.
With the Bernoulli numbers $ B_n $ defined by $ xe^x/(e^x-1) = \sum_{n=0}^\infty \frac{B_n}{n!} x^n $,
we have $ \zeta_\Q(1-n) = -\frac{B_n}n $ for all integers $ n \ge 1 $
by~\cite[Theorem~4.2]{wash}.
The theorem of Von Staudt-Clausen (see, e.g., \cite[Theorem~5.10]{wash})
gives that $ B_n + \sum_{p \textup{ prime} \atop p-1|n} \frac 1p $
is an integer when $ n \ge1 $ is even. In particular, $ |B_n|_p = p $ if $ p-1 $
divides $ n $, and $ |B_n|_p \le 1 $ otherwise.
An odd prime $ p $ is called regular if $ |B_a|_p = 1 $ for all $ a=2,4,\dots,p-3 $
and irregular if $ |B_a|_p < 1 $ for some $ a=2,4,\dots,p-3 $.

\begin{lemma}\label{vlemma}
Let $ k=\Q $, and let $ p $ be a prime number.
\begin{enumerate}
\item
$ \frac12 G_{{\bf 1}_\Q}(T) $ is in $ \Zp[[T]]^* $.
\item
Assume $ p \ne 2 $.
If $ a=2,4,\dots,p-3 $, then either
\begin{enumerate}
\item
$ |B_a|_p = 1 $, and the distinguished polynomial of $ G_{\o_p^a} $
is trivial,
or

\item
$ |B_a|_p < 1 $, and the distinguished polynomial of $ G_{\o_p^a} $
is non-trivial.

\end{enumerate}
\end{enumerate}
\end{lemma}

\begin{proof}
(1)
This is \cite[Lemma 7.12]{wash}.

(2)
If $ n $ is a positive integer with $ n \equiv a $ modulo $ p-1 $, then
$ |\frac{B_n}{n} - \frac{B_a}{a}|_p < 1 $ by~\cite[Corollary~5.14]{wash},
so $ |\frac{B_n}{n}|_p = 1 $ in (a) and $ |\frac{B_n}{n}|_p < 1 $
in~(b).
As $ H_{\o_p^a}(T) = 1 $, and
\begin{equation} \label{zetaBn}
 \Lp p {\Q} \o_p^a 1-n = \Lp p {\Q} \o_p^n 1-n = - \frac{B_n}{n} (1-p^{n-1})
\end{equation}
by~\eqref{interpol}, the two cases correspond to the function $ G_{\o_p^a}(u^s-1) $ assuming a value in $ \Zp^* $ or in
$ p\Zp $. In the first case the distinguished polynomial
of $ G_{\o_p^a}(T) $ is~$ 1 $ (and necessarily also $ \mu=0 $).
In the second case, because it is known that $ \mu=0 $ by~\cite{Fer-Wash},
the distinguished polynomial is non-trivial.
\end{proof}

From Lemma~\ref{vlemma} we see that
Conjecture~\ref{mainconjecture}(2) holds for $ k = \Q $, and
that the only \padic
\Lfunctions of the form $ \Lp p {\Q} \o_p^a s $ with $ a = 2, 4, \dots, p-3 $
for which Conjecture~\ref{mainconjecture}(3) could fail are 
those where $ p $ is an (odd) irregular prime and~$ v_p(B_a) > 0  $. 

Let us discuss the consequences of $ v_p(B_a) = 1 $, as one 
finds for an irregular prime~$ p $ in practice
(i.e., for all $ p < 12\cdot 10^6 $ as stated in
\cite[Remark~2.8]{kellner08}, or $ p < 163 \cdot 10^6 $ as stated
on \cite[p.2439]{buhler11}).
If in this case the distinguished polynomial of $ G_{\o_p^a} $ is of degree at least 2,
it is Eisenstein because $ H_{\o_p^a}(T) = 1 $,
$ v_p( G_{\o_p^a}(u^a-1) )=  v_p( \Lp p {\Q} \o_p^a 1-a ) = v_p(B_a) = 1  $,
and $ G_{\o_p^a}(u^s-1) - G_{\o_p^a}(0) $ is in $ p^2 \Zp $ for all $ s $
in~$ \Zp $ because $ u^s-1 $ is in $ p\Zp $.
Hence it has no root in $ \Zp $,
and $ \Lp p {\Q} \o_p^a s $ has no zero in $ \Zp $. In fact, then
$ | \Lp p {\Q} \o_p^a s |_p = p^{-1} $ for all $ s $ in~$ \Zp $,
so $ v_p(B_n/n) = 1 $ for all positive $ n $ congruent to~$ a $ modulo $ p-1 $
by~\eqref{zetaBn}.

But as mentioned right after~\cite[Theorem~3.2]{kellner08},
one also finds in practice (again for all $  p < 12\cdot 10^6 $)
that there is such an $ n $ with
$ v_p(B_n/n) > 1 $, so the distinguished polynomial of $ G_{\o_p^a}(T) $ must have degree~1. It has a unique (simple) root, which is in $ p\Zp $,
so that $ \Lp p {\Q} \o_p^a s $ has a unique (simple) root in $ \D_k $, which lies in~$ \Zp $
(cf.\cite[Theorem~4.10]{kellner08}).
This root is non-zero by \cite[Proposition 5.23]{wash}, so according to Conjecture~\ref{mainconjecture} it should not
be in~$ \Z $.

\begin{remark}
In Proposition~\ref{vprop} we shall see that
$ \frac12 G_{\chi}(T) $ is in $ \Zp[[T]]^* $ for any even 1-dimensional Artin character 
$ \chi : G_\Q \to \Qpbar $
of $ p $-power order with conductor~$ p^n $ for~$ n \ge 1 $,
thus generalising Lemma~\ref{vlemma}(1).
In particular, Conjecture~\ref{mainconjecture}(3) holds for such~$ \chi $.
\end{remark}

In the remainder of this section, we first discuss some earlier cases of Conjecture~\ref{mainconjecture} in
Remark~\ref{earlier}.
We then discuss the `missing' case $ m = 0 $ in Remark~\ref{zeroremark},
and conclude this section with an explanation of our motivation
for making Conjecture~\ref{Grootsconjecture}.

\begin{remark} \label{earlier}
Parts and/or special cases of Conjecture~\ref{mainconjecture} exist in the literature,
almost always for $ p $ odd. We mention a few, but this list is certainly not exhaustive.

(1)
Conjecture~\ref{mainconjecture}(2) for $ m=1 $ is equivalent
with the Leopoldt conjecture for~$ k $ (see~\cite{Colm88}).
As that is known if $ k $ is totally real and Abelian over $ \Q $
by work of Brumer~\cite{Brumer}, one sees that $ \Lp p {\Q} {\chi} 1 \ne 0 $
if $ \chi \ne {\bf 1}_\Q $ is any even Dirichlet character~\cite[Corollary~5.3]{wash}.

(2)
On \cite[p.~54]{Ba-Ne}, the hypothesis is put forward that $ \Lp p k \omega_p^i m \ne 0 $
for $ k $ totally real, $ i $ even,~$ p $~odd, $ m \ne 1 $, and $ m+i \equiv 1$ modulo the order $ d $ of $ \omega_p $.
As $ m $ is then odd, this is implied by part~(1) of Conjecture~\ref{mainconjecture}
if $ \o_p^i = {\bf 1}_{G_k} $, and by part~$ (3) $ if~$ \o_p^i \ne {\bf 1}_{G_k} $. 

(3)
In \cite[3.4]{sou81} C.~Soul\'e mentions
the thought that $ \Lp p k \omega_p^{1-m} m $ is non-zero for
odd $ m \ge 3 $ when $ k $ is a totally real Abelian extension of~$ \Q $,
and $ p $ is an odd prime not dividing $ [k : \Q] $.

(4)
For $ p $ odd,
Conjecture~\ref{mainconjecture}(2), as well as part (3) for
$ \chi = \omega_p^i  \ne {\bf 1}_{G_k} $ with $ i $ even, is considered
`reasonable' on page~996 of~\cite{mitch05}.
\end{remark}

\begin{remark} \label{zeroremark}
We discuss the `missing' case $ m=0 $ of Conjecture~\ref{mainconjecture},
where one has Conjecture~\ref{gross conj} instead.
If $ \chi $ is 1-dimensional, then from~\eqref{interpol} it is clear that~$ \Lp p k {\chi} 0 $ can
be zero, with the zeroes coming from the Euler factors on the
right-hand side as $ L^*(k , \chi\o_p^{-1} , 0) \ne 0 $, but for
a proper discussion we introduce another \Lfunction.

Let $ k $ be a totally real number field, $ i $ an integer, $ E $ a finite extension of $ \Q $, 
and~$ \psi : G_k \to E $ an Artin character with $ \psi(\sigma) = (-1)^i \psi(e) $ for
every complex conjugation~$ \sigma $ in $ G_k $.  Then one can define a
function $ \Lp p k \psi\tensor\o_p^i s $ with values in~$ E\tensor_\Q\Qpbar $
by demanding that,
for every embedding $ \tau : E \to \Qpbar $, its image under the map $ E \tensor_\Q \Qpbar \to \Qpbar $
sending $ e \tensor a $ to $ \tau(e)a $	is $ \Lp p k {\psi^\tau}\o_p^i s $ in $ \Qpbar $,
where we write $ \psi^\tau $ for $\tau \circ \psi $.
According to what is stated right after Conjecture~\ref{artinconjecture}, for $ p \ne 2 $, this function is defined
on~$ \D_k $ if $ \psi $ does not contain the trivial character, and
on~$ \D_k\setminus\{1\} $ if it does; for $ p = 2 $ it is defined on
$ \D_k $ except for finitely many points.
In fact, by (the proof of)~\cite[Lemma~3.5]{BBDJR},
this function takes values in $ E\tensor_\Q\Qp $ on the
points of $ \D_k \cap \Qp $ where it is defined.
Clearly, for an integer $ m \ne 1 $, $ \Lp p k \psi\tensor\o_p^i s $ is defined at $ s=m $ and its value
there is a unit in $ E \otimes_\Q \Qp $ if and only if, for every~$ \tau : E \to \Qpbar $,
$ \Lp p k {\psi^\tau}\o_p^i s $ is defined at $ s=m $ and its value is non-zero in $ \Qpbar $.

Now take $ i = 1 $.
Using all the embeddings of
$ E $ into~$ \Qpbar $ and \eqref{iwasawaseries}, one sees that
there is a Laurent series expansion (and in fact, a power series
expansion if $ p \ne 2 $) around~0 in $ \D_k $ with coefficients
in $ E \otimes \Qp $.

Let $ r $ be the order of vanishing at $ s = 0 $ of the truncated complex $L$-function
\begin{equation*}
  L(k , \psi^\tau , s) \prod_{\mathfrak{p} \in P}  \Eul_\mathfrak{p}(k, \psi^\tau,s) 
\,,
\end{equation*}
for any embedding $ \tau : E \to \C $. As mentioned above, the vanishing of this function at $s=0$ comes uniquely from the Euler factors.
More precisely, this is determined by when a certain decomposition group is contained in the kernel of $\psi^\tau$; see the discussion before Theorem~\ref{th:proofgoss=0}. Since this kernel is independent of $\tau$, the order of vanishing $r$ is also independent of~$\tau$.  
Let $ U(s) s^{r'} $ with $ U(s) $ in $ E\otimes_\Q\Qp[[s]] $ and $ U(0) \ne 0 $ be the expansion around~0 of the $p$-adic $L$-function $ \Lp p k \psi\tensor\o_p s $; thus~$r'$ is the order of vanishing at $ s = 0 $ of $ \Lp p k \psi\tensor\o_p s $.
Gross conjectured (see \cite[Conjecture~2.12]{Gross81}) that~$r' = r$, and gave a precise formula for $ U(0) $, which should be in $ (E\otimes_\Q\Qp)^* $. It follows from the interpolation property of $p$-adic $L$-functions that $r' = r$ if $r= 0$. Equality
also holds when $r = 1$ by a result of Gross \cite[Proposition~2.3]{Gross81}. 
Thanks to this result, we proved numerically that Conjecture~\ref{gross conj} holds in many cases; see Theorem~\ref{th:proofgoss=0}. 
For the general case, Spie\ss~\cite{Spiess14} proved that $r' \geq r$ when $\psi$ is $1$-dimensional and, recently, the conjectural formula for $\tfrac{1}{r!} \Lp p\egroup^\bgroup(r) k \psi\tensor\o_p 0 $
was proved for $ \psi $ 1-dimensional in~\cite{DaKaVe}. However, it is not known that this value is non-zero
in general
for $ 1 $-dimensional $ \psi $, and therefore the inequality $r' \leq r$ is still open. 
In these cases, to check that $ r = r' $ and $ U(0) $ is in $ (E \otimes_\Q \Qp)^* $, we just need to compute $\Lp p\egroup^\bgroup(r) k \psi\tensor\o_p 0 $ with enough precision to ensure that it is in $ (E\otimes_\Q \Qp)^* $. We refer to Remark~\ref{grossremark} for an example of such a verification of the conjecture with $r = 2$. 

\end{remark}

We conclude this section with a discussion of Conjecture~\ref{Grootsconjecture}.
This conjecture is based on the complex
situation. Indeed, assume that $\psi : G_k \to \C$ is a 1-dimensional
Artin character of $k$ such that $\Ind {G_k} {G_\Q} (\psi)$
is irreducible. A consequence of the (classical) Artin conjecture
is that the complex $L$-function $L(\Q, \Ind {G_k} {G_\Q} (\psi), s) = L(k, \psi, s)$ should
be a primitive element of the Selberg class.
Furthermore, one expects that a primitive function in the
Selberg class has at most one non-trivial zero of multiplicity $> 1$,
occurring at the point $s=1/2$; see for example the discussion after
Conjecture~7.1 in \cite[p.\ 107]{per04}. The point $s=1/2$ plays a
special role in the complex case because of the functional equation,
but it does not have an equivalent in the $p$-adic world. Thus, the
natural equivalent in the $p$-adic world is to conjecture that the
$p$-adic $L$-function $\Lp p k {\chi} {s} $ has only simple
zeroes except for the possible trivial zero at~$s=0$; see Remark~\ref{zeroremark}.
In Conjecture~\ref{Grootsconjecture} we used the stronger statement that $D_\chi(T)$ has only simple zeroes apart from the possible trivial zero at~$T=u-1$.

\begin{remark}
In \cite[Conjecture~5.19]{dJNa}(1) a related but different conjecture is made
for an even 1-dimensional Artin character $ G_k \to \Qpbar $
of order prime to~$ p $.
It states that in that case the distinguished polynomials $ g_j(T) $ among
the invariants of the corresponding Iwasawa module are square-free.
If $ p \ne 2 $, or if $ p = 2 $ under Assumption~4.4 of loc.~cit.,
we have $ G_\chi(T) = \prod_j g_j(T) $.
Thus, if also Conjecture~\ref{Grootsconjecture} applies, then
the $ g_j(T) $ can only have a factor $ T-u+1 $
in common, and this factor can occur at most once in every
$ g_j(T) $.
\end{remark}

\begin{remark}
Note that the discussion following Lemma~\ref{vlemma} shows that
there is much numerical evidence for both Conjecture~\ref{Grootsconjecture}
and \cite[Conjecture~5.19]{dJNa} in the case that $ k = \Q $,
$ p $ is odd, and $ \chi = \o_p^a $ for $ a = 2, 4, \dots, p-3 $.
\end{remark}

%% file: consequences.tex
\section{Some consequences of Conjecture~\ref{mainconjecture} for \etale cohomology and for syntomic regulators} \label{consequences}

In this section we discuss some very important consequences of Conjecture~\ref{mainconjecture},
for \etale cohomology groups through Theorems 1.4 and 1.8 of~\cite{dJNa},
and for $ p $-adic (syntomic) regulators through Conjecture~3.18
and Proposition~4.17 of~\cite{BBDJR}.
The latter consequences were in fact the motivation for making
Conjecture~\ref{mainconjecture}.

In order to state the consequences for the \etale cohomology groups, we quote Theorem~1.8(3) of~\cite{dJNa}
as Theorem~\ref{dJNathm} below.
In fact, loc.\ cit.\ describes a more general situation,
using a modified Tate twist with a \padic index $ 1-e $
instead of an integral $ 1-m $, but we do not need this kind of generality here.

For the notation used in Theorem~\ref{dJNathm}, we first recall
that the truncated \padic \Lfunction $ L_{p,S} $ was discussed at
the end of Section~\ref{Lpsection}.
For the remainder of the theorem, suppose $ k $ is a (not necessarily totally real) number field, $ E $ a finite extension of $ \Qp $,
and $ \eta : G_k \to E $ the Artin character 
of some representation of~$ G_k $ over $ E $.
With $ \Oe $ the valuation ring of $ E $, one can find a 
finitely generated, torsion-free
$\Oe$-module $M_E(\eta) $ on which $G_k$ acts and such that the
resulting representation of $ G_k $ on  $M_E(\eta) \otimes_{\Oe} E \simeq M_E(\eta) \otimes_{\Zp} \Qp $
has character~$ \eta $.
We then put~$ W_E(\eta) = M_E(\eta) \otimes_{\Zp} \Qp/\Zp $.

\begin{theorem} \label{dJNathm}
Let $ k $ be totally real, $ p $ a prime, $ E $ a finite extension
of $ \Qp $, and assume the Artin character $ \chi $
of $ G_k $ comes from a representation of $ G_k $
over~$ E $ and is even.
Let $ S $ be a finite set of finite primes of $ k $ containing
the primes above $ p $,
as well as the finite primes at which $ \chi $ is ramified, and
let $ \O_{k,S} $ be obtained from the ring of
integers $ \O_k $ of~$ k $ by inverting all primes in~$ S $.
Let $ m $ be an integer, with $ m \ne 1 $ if $ \chi $ contains the trivial character $ {\bf 1}_{G_k} $.
Finally, for $ p=2 $, assume the slightly weaker version
of the main conjecture of Iwasawa theory as formulated in \cite[Assumption~4.4]{dJNa}.
Then the following are equivalent.
\begin{enumerate}
\item[(i)]
$ \Lp p,S k {\chi} m \ne 0 $;

\item[(ii)]
$H_{\et}^1(\O_{k,S} , W_E(\chi\o_p^{m-1})(1-m)) $ is finite;

\item[(iii)]
$H_{\et}^2(\O_{k,S} , W_E(\chi\o_p^{m-1})(1-m)) $ is finite;

\item[(iv)]
$H_{\et}^2(\O_{k,S} , W_E(\chi\o_p^{m-1})(1-m)) $ is trivial.
\end{enumerate}
If these equivalent conditions hold then
\begin{equation} \label{multeuler}
| \Lp p,S k {\chi}  m |_p ^{[E:\Q_p]} = 
\frac{ \# H_{\et}^0(\O_{k,S},W_E(\chi\o_p^{m-1})(1-m)) }{ \# H_{et}^1(\O_{k,S},W_E(\chi\o_p^{m-1})(1-m) ) }
\,. 
\end{equation}
\end{theorem}

For $ m \ne 0 $, part (i) of this theorem is equivalent to $ \Lp p k {\chi} m \ne 0 $
because the Euler factors in~\eqref{pLS} are non-zero for $ m \ne 0 $ as mentioned
at the end of Section~\ref{Lpsection}.
So if $ m < 0 $ then (i) holds by~\eqref{interpol}, and~\cite[p.2338-2339]{dJNa}
then used~\eqref{multeuler} and~\eqref{pLS} to 
prove a (generalisation of a) conjecture
by Coates and Lichtenbaum \cite{Co-Li} for $ p $ odd.

If $ k $ is fixed, and $ m \ge 2 $ (or $ m \ge 1 $ if $ \chi $ does not contain $ {\bf 1}_{G_k} $), then
(i) above is implied by part~(4) (or parts~(2) and~(4)) of Conjecture~\ref{mainconjecture}
for the same $ k $ and $ m $. As a consequence of Theorem~\ref{dJNathm}, then the conjecture of Coates
and Lichtenbaum in its incarnation~\eqref{multeuler} would hold for such $ k $, $ m $ and $ \chi $,
and odd~$ p $.

The numerical calculations we performed prove that $ \Lp p k {\chi} m \ne 0 $
for many~$ p $,~$ k $ and~$ \chi $, with $ k $ equal to $ \Q $ or a real quadratic field,
either for all $ m \ne 0 $ or for all $ m \ne 0 $ up to a large positive
bound (see Theorem~\ref{th:mainconj} for the precise statement).
In particular, in those cases the above equivalent statements
(i)-(iv) and~\eqref{multeuler} hold.

Perhaps superfluously, we note that the three previous paragraphs also apply for~$ p = 2 $
if we make~\cite[Assumption~4.4]{dJNa}.

\medskip

We now discuss the consequences of Conjecture~\ref{mainconjecture}
for a conjecture by Schneider.
He conjectured (see \cite[p.~192]{schn79})
that, for $ n \ne 1 $ an integer, and $ k $ any number field,
the maximal divisible subgroup of $ H_{\et}^2(\O_k[p^{-1}] , \Qp/\Zp(n)) $
is trivial if $ p $ is an odd prime, but for simplicity we discuss
the case $ p = 2 $ as well. By \cite[Theorem~1.4(4)]{dJNa} this is equivalent
with this \etale cohomology group being finite (and, in fact,
trivial if $ p \ne 2 $).

For $ k $ totally real and $ n \ne 0 $ even, and making the assumption
in Theorem~\ref{dJNathm} if $ p = 2 $, by that theorem 
(with $ E = \Qp $, $ M_E = \Zp $, $ \chi = \o_p^n $ and~$ m = 1-n $),
this is equivalent to $ \Lp p k {\o_p^n} 1-n \ne 0 $, which is implied by parts~(2)
and~(3) of Conjecture~\ref{mainconjecture}.
For $ k $ totally real and $ n = 0 $, under the same assumption
if~$ p = 2 $, using \cite[Remark~5.20(1)]{dJNa} one sees that
it is equivalent to $ \zeta_p(k,s) $ having a simple pole at
$ s = 1 $, i.e., to part of Conjecture~\ref{mainconjecture}(2).
We refer to Remark~\ref{rk:schneider} for our results on Schneider's conjecture.

\medskip

In order to discuss the consequences of Conjecture~\ref{mainconjecture} for syntomic regulators, we
use the \Lfunctions introduced in Remark~\ref{zeroremark}, but
only for $ k = \Q $. We let $ E $ be a finite extension of $ \Q $, $ m \ge 2 $
an integer, and~$ \psi : G_{\Q} \to E $ an Artin character with~$ \psi(\sigma) = (-1)^{1-m} \psi(e) $ for
every complex conjugation~$ \sigma $ in $ G_{\Q} $.

According to \cite[Conjecture~3.18(4)]{BBDJR}, $ \Lp p {\Q} \psi\tensor\o_p^{1-m} m $ 
should exist and be a unit in $ E \tensor_\Q \Qp $.
The condition on $ \psi $ means that $ \psi^\tau \o_p^{1-m} $ is even for one (and hence
every) embedding
$ \tau : E \to \Qpbar $, so that the statement
is equivalent to~$ \Lp p {\Q} {\psi^\tau\o_p^{1-m}} m $ existing and being non-zero
for all $ m \ge 2 $ and even $ \psi^\tau \o_p^{1-m} $.
Since every even Artin character $ \chi : G_\Q \to \Qpbar  $ is of this form,
and $ \zeta_{\Q,p}(n) \ne 0 $ for every integer $ n \ne 1 $ by
Lemma~\ref{vlemma}(1), the statement of loc.\ cit.\ is equivalent to Conjecture~\ref{mainconjecture}(5)
for~$ m \ge 2 $.

In fact, \cite[Conjecture~3.18(4)]{BBDJR} was one of the main reasons for making Conjecture~\ref{mainconjecture},
as together with part~(2) of that conjecture, it implies that
certain~\padic\ (syntomic) regulators should be units in (a subalgebra of)
$ E \otimes_\Q \Q_p $, in analogy to the Leopoldt conjecture (cf.\ Remark~\ref{earlier}(1)
in the current paper).

According to \cite[Proposition~4.17]{BBDJR}), parts~(1) through~(3)
of Conjecture~3.18 of loc.~cit.~hold if $ \psi $ is a 1-dimensional character of $ G_\Q $ for
which the parity of $ \psi $ and of~$ 1-m $ match. 
Our calculations prove that part~(4) of
this conjecture then also holds for all $ p $, $ m \ge 2 $ and $ \psi $ such that all
$ \chi = \psi^\tau \o_p^{1-m} $ occur in  Theorem~\ref{th:mainconj}(1).

Moreover, for certain even $ \chi : G_\Q \to \Qpbar $ it follows from the theory that
$ \Lp p k {\chi} s $ does not have any zero (see Lemma~\ref{vlemma}
and Proposition~\ref{vprop}.)
Hence Conjecture~3.18 of loc.~cit.~holds for the corresponding
$ \psi $ in full in those cases,
and in particular, all corresponding regulators are units in
the corresponding $ E \otimes_\Q \Qp $.
We refer to Corollary~\ref{regcor}
for a precise description of those $ \psi $.

%% file: NP.tex
\section{Power series and Newton polygons} \label{NPsection}

In this section we collect some information about power series
and Newton polygons, which will be used in the description of
the calculations in Section~\ref{zeronumevidence}.

Let $ E $ be a finite extension of $ \Qp $, with valuation
ring $ \O $ and residue field~$ \kappa $.
Let~$ S(T) = \sum_{i \ge 0} S_i T^i $ in $ \O[[T]] $ be such that its reduction $ \overline{S(T)} $
in $ \kappa[[T]] $ is non-zero, and let $ d \ge 0 $ be the lowest
index such that $ S_d $ is a unit of~$ \O $. 
As stated just after~\eqref{iwasawaseries}, we can write 
$ S(T) = D(T) U(T) $
with $ D(T) $ a distinguished polynomial of degree $ d $, and $ U(T) $
in $ \O[[T]]^* $.
Moreover, writing $ S(T) $ as a product of a distinguished polynomial
and a unit of $ \O[[T]] $ can be done in only this way, so,
in particular, they remain the same if we enlarge~$ E $.

\medskip

 We now briefly recall the definition and first properties of Newton
polygons following Koblitz \cite[IV.3-4]{Kob77}. Note that Koblitz
assumes that the power series under consideration all have constant
term equal to~1, whereas we will only assume that they have non-zero
constant term.

\begin{figure}[htp]
\begin{tikzpicture}[scale=0.65]
\draw (0,7) -- (0,0) -- (17,0);
\foreach \Point in {(0,6), (3,5), (4,4), (5,5), (6,3), (8,4), (11,2), (12,1), (15,3), (16,0)}{ \node at \Point {\textbullet};}; 
\draw (0,6) -- (6,3) -- (12,1) -- (16,0); 
\end{tikzpicture}
\caption{$\mathcal{NP}(S)$ with $S(T) = p^6 + p^5 T^3 + p^4 T^4 + p^5 T^5 + p^3 T^6 + p^4 T^8 + p^2 T^{11} + p T^{12} + p^3 T^{15} + T^{16} + \cdots$. The polynomial $P(T)$ associated to the first segment is $1 + T^4 + T^6$.}\label{fig:2}
\end{figure}
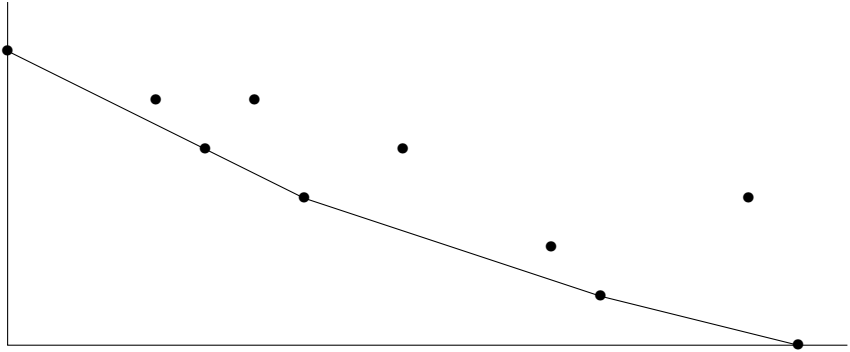

For this, we now assume that $S_0 \ne 0$. We refer to \cite{Kob77}
for the exact definition of the Newton polygon of $S$; loosely
speaking it is the ``convex hull'' of the points~$(i, v_p(S_i))$,
for $i \geq 0$, and is obtained by starting with a vertical half-line
through $(0, v_p(S_0))$ and rotating it about this point
counterclockwise until it hits one of the points
$(i_1, v_p(S_{i_1}))$ with $i_1 > 0 $ and $i_1$ maximal if several points are hit, giving the first segment joining
$(0, v_p(S_0))$ and $(i_1, v_p(S_{i_1}))$,
then rotating the remainder of the half-line
further about $(i_1, v_p(S_{i_1}))$ until it hits another point~$(i_2, v_p(S_{i_2}))$ with $i_2 > i_1$ and, as above, $i_2$ chosen maximal if several values are possible, giving the second segment
joining $(i_1, v_p(S_{i_1}))$ and $(i_2, v_p(S_{i_2}))$, etc.\ (see
Figure~\ref{fig:2}). 
This process either goes on indefinitely, or terminates with a segment of infinite length such that it can be rotated no further without going beyond at least one point~$ (i, v_p(S_i)) $ but not hitting any of those points.

We need some
additional definitions and notations. Consider a segment $\sigma$ of
the Newton polygon joining the points $(i, y)$ and $(j, z)$ with $i < j$.
We call $(i,y)$ the origin of the segment and $(j,z)$ the end of
the segment; the height $h(\sigma)$ of $\sigma$ is defined by
$h(\sigma) = z-y$; the length $l(\sigma)$ is defined by $l(\sigma) = j-i$;
the slope $s(\sigma)$ of $\sigma$ is defined by $s(\sigma) = h(\sigma)/l(\sigma)$.

All the power series and polynomials that we
consider have coefficients with non-negative valuations, non-zero constant term, and at
least one coefficient with valuation zero.
We use the convention that the Newton polygon of
a power series, or if it applies, the Newton polygon of a polynomial,
ends with the segment that reaches the horizontal axis. Therefore the
Newton polygons we consider all have finitely many segments, all with
negative heights and slopes, and the slopes are strictly increasing.
We will denote by $\mathcal{NP}(S)$ the Newton
polygon of $S$, which we see as the set of its segments.

Assuming that $ S(T) $ is in $ \O[[T]] $, has non-trivial constant
term, and does not reduce to zero in $ \kappa[[T]] $,
we want to identify the Newton polygons of~$ S(T) $ and~$ D(T) $,
together with some additional information (see Lemma~\ref{linklemma}
below). For this we introduce the following.

\begin{definition}
We shall call $ (i, v_p(S_i)) $ a  point on the Newton polygon of
$ S(T) $ if it lies in one of the segments in $ \mathcal{NP}(S) $.
\end{definition}

Note that the origin and end of each segment are points on the
Newton polygon.
Clearly, if $ D(T) = 1 $ then $ \mathcal{NP}(D) $ and $ \mathcal{NP}(S) $ are
both empty,
so assume that~$ D(T) $ is a distinguished polynomial of positive degree $ \l $.
Any irreducible factor in~$ \O[T] $ will have roots of the same
positive valuation $ -s $ for some negative rational number~$ s $.
Collecting factors with the same valuation
for its roots together, we can write
\begin{equation*}
 D(T) = D_{s_1}(T) \cdots D_{s_l}(T) 
\end{equation*}
where $ s_1 < \dots < s_l $, each $ s_j $ is a negative rational
number, and $ D_{s_j}(T) $ is of positive degree $ \l_j $. Note
that each $ D_{s_j}(T) $ is again distinguished, and that
$ U(T) $ has no roots of positive valuation, so such roots of
$ S(T) $ are precisely the roots of the~$ D_{s_i}(T) $.

It is easy to verify that the Newton polygon of each $ D_{s_j}(T) $
is a segment of length~$ \l_j $, and height $ s_j \l_j $, with
vertices $ (0,-s_j \l_j) $ and $ (\l_j,0) $.
Now assume that the Newton polygon of $ D_{\le s_j}(T) := D_{s_1}(T) \cdots D_{s_j}(T) $
is known for some $ j $ with $ 1 \le j \le l-1 $,
and that it consists of $ j $ segments  with lengths $ \l_1,\dots,\l_j $
and slopes $ s_1,\dots,s_j $. We want to
obtain from it the Newton polygon of $ D_{\le s_{j+1}}(T) =  D_{\le s_j}(T) D_{s_{j+1}}(T) $.
Note that the Newton polygon of $ D_{\le s_j}(T) D_{s_{j+1}}(0) $
is just the vertical translate of the one for $ D_{\le s_j}(T) $ over
$ v_p(D_{s_{j+1}}(0)) $.
Similarly, the term $ T^{\l_1+\dots+\l_j} D_{s_{j+1}}(T) $ gives
the translate of the Newton polygon of $ D_{s_{j+1}}(T) $ to
the right over $ \l_1+\dots+\l_j $. Both translates meet at the
point $ (\l_1+\dots+\l_j, v_p(D_{s_{j+1}}(0)) ) $. Because $ s_{j+1} $
is larger than the slopes occurring in the Newton polygon for
$ D_{\le s_j}(T) $, it is easy to check that other terms are
above the concatenation of the two translates. Therefore,
this concatenation is the Newton polygon of $ D_{\le s_{j+1}}(T) $,
so the only slopes that occur are now $ s_1,\dots,s_{j+1} $,
with lengths $ \l_1,\dots, \l_{j+1} $.
Moreover, this calculation also shows that the only points on the Newton polygon
are the ones that come from $ D_{\le s_j}(T) D_{s_{j+1}}(0) $
and  $ T^{\l_1+\dots+\l_j} D_{s_{j+1}}(T) $, i.e., they are translates
of points on the Newton polygons from $ D_{\le s_j}(T) $ and
$ D_{s_{j+1}}(T) $. By induction, we see that the only slopes
that occur in the Newton polygon of $ D(T) $ are $ s_1,\dots,s_l $,
with lengths $ \l_1,\dots,\l_l $, and that the only points on
it are the ones corresponding (under the translations) to the
ones on the Newton polygons for the $ D_{s_j}(T) $.

Multiplying $ D(T) $ by $ U(T) $ in $ \O[[T]]^* $
does not change the Newton polygon under our conventions, because $ D(T) U(0) $ has the same
Newton polygon as $ D(T) $, and multiplication by anything in
$ T\O[[T]] $ will gives us points to the right or above the Newton
polygon because all slopes in it are negative. This also shows
that the points on the Newton polygons for $ D(T) $ and $ S(T) = D(T) U(T) $
are the same.

The following lemma will enable us to rule out multiple roots
of $ S(T) $ based on the points on $ \mathcal{NP}(S) $
because the roots of $ S(T) $ in $ \Qpbar^* $ with positive valuation
are exactly those of $ D(T) $.

\begin{lemma} \label{linklemma}
Let $ S(T) $ in $ \O[[T]] $ be such that it has non-trivial constant
term, and does not reduce to zero in $ \kappa[[T]] $.
Fix a uniformizing parameter $ \pi $ of $\O$ and let~$ e = v_p(\pi)^{-1} $
be the ramification index of~$ E $.
Let $ \sigma = \sigma_i $ be one of the segments of the Newton polygon
of $ S(T) $, with length $ \l = \l_i $ and slope~$ s = s_i $.
If $ \sigma $ has origin $ (c,v_p(S_c)) $ and end $ (c+\l, v_p(S_{c+\l})) $,
Let $ \sum_{j=c}^{c+\l} S_j' T^j $
be the sum of the terms of $ S(T) $ corresponding to the points on~$ \sigma $,
i.e.,
\begin{equation*}
 S_j' = \left\{
\begin{aligned}
   S_j \text{ if } v_p(S_j)-v_p(S_c) = (j-c) s \\
   0 \text{ if } v_p(S_j)-v_p(S_c) > (j-c) s
\,.
\end{aligned}
\right.
\end{equation*}
Consider the polynomial 
\begin{equation*}
 P_\sigma(T) = T^{-c} \sum_{j=c}^{c+\l} \overline{S_j' \pi^{-ev_p(S_j')}} T^j 
\end{equation*}
in $ \kappa[T] $.
If $ P_\sigma(T) $ has no multiple factors in $ \kappa[T] $, then $ D_s(T) $
has no multiple factors, so that, in particular, $ S(T) $ has no multiple
roots in $ \Qpbar^* $~of slope $ s $.
\end{lemma}

\begin{proof}
Write $ -e s = a/b $ for positive integers $ a $ and $ b $.
Fix some $ \tilde \pi $ in $ \Qpbar $ with~$ \tilde \pi^b = \pi $.
Then $ \tilde E = E(\tilde\pi) $ is purely ramified of degree
$ b $ over~$ E $, and $ \tilde\pi $ is a uniformizing parameter
for its valuation ring~$ \tilde\O $.
Let $ \rho = \tilde \pi^a $, so $ v_p(\rho) = \frac ab v_p(\pi) = -s > 0 $.

Write $ D_{<s}(T) = D_{s_1}(T) \dots D_{s_{i-1}}(T) $
and $ D_{>s}(T) = D_{s_{i+1}}(T) \dots D_{s_l}(T) $, so that
we have $ S(T) = D_{<s}(T) D_s(T) D_{>s}(T) U(T) $ with $ U(T) $ in $ \O[[T]]^* $
by our assumptions.
So $ \l = \deg(D_s) $, $ c = \deg(D_{<s}) $, the origin of
$ \sigma $ is $ (c , v_p( D_s(0) D_{>s}(0)) ) $ and its end $ ( c + \l , v_p(D_{>s}(0)) ) $.

Let $ A(T) = \sum_{n=0}^\infty A_n T^n =  \rho^{-c-\l} D_{>s}(0)^{-1} S(\rho T) $.
From the valuations of the coefficients of $ S(T) $ one sees
that $ A(T) $ is in $ \tilde \O[[T]] $, and that for $ A_n \ne 0 $ we have~$ v_p(A_n) = 0 $ if and only if $ (n, v_P(S_n)) $ is a point
on the segment $ \sigma $ of $ \mathcal{NP}(S) $.
Hence reducing its coefficients modulo $ \tilde \pi \tilde \O $ gives
the same as reducing the coefficients of $ \rho^{-c-\l} D_{>s}(0)^{-1} \sum_{j=c}^{c+\l} S_j' (\rho T)^j $.
For those $ j = c,\dots, c+\l $ for which $ S_j' \ne 0 $ we have $ v_p(S_j') - v_p(S_{c+\l}') = s (j - c - \l) $,
so that
\begin{equation*}
 \rho^{-c-\l} S_j' \rho^j 
=
  S_j' \tilde \pi^{a (j-c-\l)}
=
 S_j' \pi^{-e s (j-c-\l)}
=
 S_j' \pi^{-e v_p(S_j')} \pi^{e v_p(S_{c+\l}')}
\,.
\end{equation*}
From this we see that reducing the coefficients of $ A(T) $ modulo
$ \tilde \pi \tilde \O $ gives $ u T^c P_\sigma(T) $,
with 
$ u $ in $ \kappa^* $ the reduction modulo $ \tilde \pi \tilde \O $
of $ D_{>s}(0)^{-1} \pi^{e v_p(S_{c+\l}')} $ in $ \tilde \O^* $.

On the other hand, we can write
\begin{equation*}
A(T)
 =
 [\rho^{-c} D_{<s}(\rho T)] [\rho^{-\l} D_s(\rho T)] [D_{>s}(0)^{-1} D(\rho T)] U(\rho T)
\,.
\end{equation*}
Because the roots of $ D_{<s}(T) $ all have valuation larger
than $ -s $, the first term here lies in $ T^c + \tilde \pi \tilde \O [T] $.
The third term is in $ 1 + \tilde \pi \O [T] $ because the
roots of $ D_{>s} $ all have valuation less than $ -s $.
The fourth term is in $ \tilde \O [[T]]^* $ as~$ v_p(\rho) > 0 $,
and reducing its coefficients modulo~$ \tilde \pi \tilde \O $
gives the same as $ U(0) $.
The second term is a monic polynomial~$ P_s(T) $ in~$ \tilde \O [T] $ with constant
term in~$ \tilde \O^* $ because all its roots have valuation~0.
So reducing the coefficients of $ A(T) $ modulo $ \tilde \pi \tilde \O $
results in $ T^c \ol P_s(T) \ol{U(0)} $, with~$ \ol P_s (T) $
obtained this way from $ P_s (T) $.

So, up to multiplication by an element of $ \kappa^* $,
we have $ P_\sigma(T) = \ol  P_s (T) $. Our claim follows because
associating $ \ol P_s(T) $ to $ D_s(T) $ is compatible with
factorising~$ D_s(T) $.
\end{proof}

\begin{figure}[h]
\begin{tikzpicture}[scale=0.95]
\draw (0,5) -- (0,0) -- (12,0);
\foreach \Point in {(0,4), (2,3), (3,2), (6,2), (9, 1)}{ \node at \Point {\textbullet};}; 
\draw (0,4) -- (3,2);
\draw[dotted] (3,2) -- (10,0);
\draw (10,-0.2) -- (10,2.5) node[anchor=west] {$L_1 = 10$} -- (10,4.5);
\end{tikzpicture}

\caption{$\mathcal{NP}(\tilde{S})$ for $\tilde{S}(T) = p^4 + p^3T^2 + p^2T^3 + p^2 T^6 + p T^9$.
Assume that $S(T) \equiv \tilde{S}(T) \pmod{T^{10}}$, then a priori only the first segment of $\mathcal{NP}(\tilde{S})$ belongs to $\mathcal{NP}(S)$ and the remaining segments of $\mathcal{NP}(S)$ have slope $\ge -2/7$.}\label{fig:np}
\end{figure}
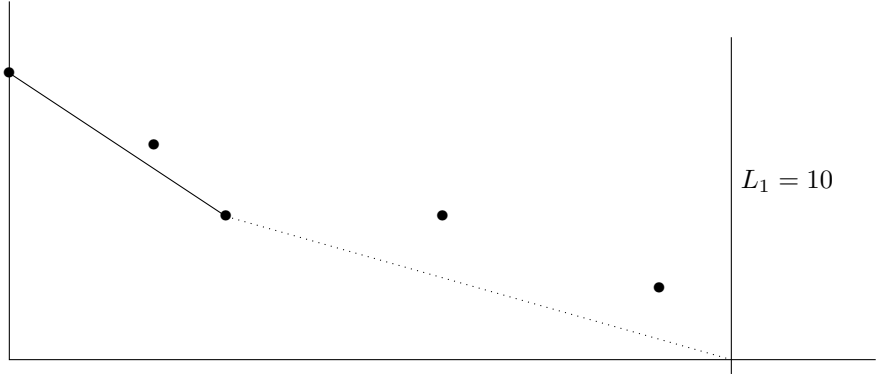

We will have to deal with the case where we only know an approximation, say $\tilde{S}(T)$, of the power series $S(T)$ modulo $T^{L_1}$ and we want to deduce information on $\mathcal{NP}(S)$ from $\mathcal{NP}(\tilde{S})$. For that, we use the following lemma whose proof is immediate (see also Figure~\ref{fig:np} for an example).

\begin{lemma}\label{lem:incnp}
Let $L_1 \geq 1$ be an integer and assume that $\tilde{S}(T) \in \O[[T]]$
is such that $\tilde{S}(T) \equiv S(T) \pmod{T^{L_1}}$ and that $S(T)$ has non-trivial constant term and does not reduce to zero in $\kappa[[T]]$. Denote by $\sigma_1, \sigma_2, \dots$ the ordered segments of $\mathcal{NP}(\tilde{S})$. Let $n \geq 1$ be the largest index such that the origin $(i,y)$ of the segment $\sigma_n$ satisfies $i < L_1$ and its slope $s$ satisfies $-y/(L_1-i) > s$.
Then the segments $\sigma_1, \dots, \sigma_n$ belong to $\mathcal{NP}(S)$. Let $(j,z)$ be the end of $\sigma_n$. If $z = 0$ then we fully obtain $\mathcal{NP}(S)$ in that way. Otherwise, the remaining segments in $\mathcal{NP}(S)$ have slope $\geq -z/(L_1-j)$. \qed
\end{lemma}

\begin{remark}
There is a segment of slope $-z/(L_1-j)$ in the Newton polygon of $S(T)$ if and only if the coefficient of $T^{L_1}$ in $S(T)$ is a unit and then it is the last slope of the Newton polygon. Otherwise, all the remaining segments in $\mathcal{NP}(S)$ have strictly larger slope. 
\end{remark}

%% file: zeronumevidence.tex
\section{Numerically verification of several conjectures on \padic \Lfunctions} \label{zeronumevidence}

In this section, we explain how we computed numerically the zeroes of the $p$-adic $L$-function associated to a $ 1 $-dimensional even $ p $-adic Artin character, for many such characters, in order to prove Conjecture~\ref{mainconjecture}(3)
for many $ k = \Q $ or a real quadratic field, $ \chi \ne {\bf 1}_k $, and $ m \ge 1 $.

Previous works have been devoted to the computation of zeroes of $p$-adic Artin $L$-functions; see for example the extensive computations in \cite{Er-Me92} and \cite{Er-Me94} and the bibliography quoted in these papers. Similarly to these previous works, we deduce
approximations of these zeroes from computed approximations of the power series $G_\chi(T)$. However, one important difference is that they computed an approximation of $ G_\chi(T) $ from an approximation of the Taylor expansion at $s=0$ of the $p$-adic Artin $L$-function, whereas we compute it directly. In particular, we  compute good approximations of the power series $G_\chi(T)$ and use them to prove Conjectures~\ref{Grootsconjecture} for those $ \chi $ as well.
We also use them to provide evidence for Conjecture~\ref{muconjecture}.

For our computations, we use the computer system~PARI/GP \cite{PARI2}. Many tools for computations in class field theory are available in this system. Thus it makes more sense to work with Hecke characters than with $1$-dimensional Artin characters; of course, these are essentially the same thing as we now explain.
 
Let $k$ be a number field. Recall that a Hecke character $ \hchi $ of $k$ is a (multiplicative) character on some ray class group  $\mathrm{Cl}_k(\mathfrak{m})$ of $k$ of modulus $\mathfrak{m}$. The character is primitive if it cannot be factored through a ray class group  $\mathrm{Cl}_k(\mathfrak{n})$ with $\mathfrak{n} \mid \mathfrak{m}$ and $\mathfrak{n} \not= \mathfrak{m}$; the modulus $\mathfrak{m}$ is then the conductor of $\hchi$.

Let $\chi : G_k \to \Qbar$ be a  $1$-dimensional Artin character. Let $\mathfrak{m}$ be the conductor of $\chi$, that is the smallest modulus (for division) such that $k_\chi \subset k(\mathfrak{m})$, where $k(\mathfrak{m})$ denotes the ray class field of $\mathfrak{m}$. Denote by $\mathrm{Art}_\mathfrak{m} : \mathrm{Cl}_k(\mathfrak{m}) \to \Gal(k(\mathfrak{m})/k)$ the Artin map,
which is an isomorphism in this case. 
Then $\hchi = \chi \circ \mathrm{Art}_\mathfrak{m}$ is a primitive Hecke character of $k$ of conductor $\mathfrak{m}$. We say that the Artin character $\chi$ and the Hecke character $\hchi$ are associated. It is immediately
seen that this construction, for all $ \mathfrak{m} $, yields a $1$-to-$1$ correspondence between the set of $1$-dimensional Artin characters of $G_k$ of conductor $\mathfrak{m}$  and the set of primitive Hecke characters of $k$ of conductor $\mathfrak{m}$. 

Assume now that $k$ is totally real. Let $\chi$ be a $1$-dimensional $ \Qpbar $-valued Artin character of $G_k$,
let $\hchi$ be the associated Hecke character, and let $\mathfrak{m}$ be
their conductor. Then $\chi$ is even if and only if $\mathfrak{m}$ has no infinite part, that is $\mathfrak{m} = \mathfrak{f}$ where $\mathfrak{f}$ is an integral ideal.
Moreover, we have that $\Lp p k {\chi} s = \Lp p k {\hchi} s $,
where $\Lp p k {\hchi} s $ is the $p$-adic Hecke $L$-function of $\hat\chi$ constructed in  \cite{Cas-Nog79} or \cite{Del-Rib}. 
For many such primitive $p$-adic Hecke characters $\hchi$ defined over $\Q$
(that is $p$-adic Dirichlet characters) or over some real quadratic field,
we compute an approximation of the power series~$G_\chi(T)$ for
the associated Artin character~$ \chi $.

We now explain which Hecke characters $\hchi$ we considere and how we construct them. The computation of the $p$-adic $L$-functions $\Lp p k {\hchi} s $ is done using the method of \cite{Rob15}.  Although the method works for any totally real number field $ k $, it is at the moment practical only when the degree of $k$ is at most $2$. Therefore we restrict to $k = \Q$ or $k$ a real quadratic field. Furthermore, if $k$ is a real quadratic field, we consider only characters $\hchi$ such that the field $k_{\hchi}$ corresponding to the kernel of $\hchi$ by class field theory, thus $k_{\hchi} = k_\chi$, is not Abelian over $\Q$. Indeed, if~$k_{\hchi}/\Q$ is Abelian then the $p$-adic $L$-function of $\hchi$ is the product of distinct $p$-adic $L$-functions of characters of $\Q$ and thus the validity for $\hchi$ of Conjecture \ref{mainconjecture}, \ref{gross conj}, and \ref{muconjecture}, reduces to the case $k = \Q$.  Note that Conjecture~\ref{Grootsconjecture} does not apply in this case since $\Ind {G_k} {G_\Q} (\chi)$ is reducible. 
We can identify the characters $\hchi$ such that $k_{\hchi}/\Q$ is Abelian thanks to the following lemma.

\begin{lemma}
Assume $k$ is a quadratic field. Denote by $\tau$ the non-trivial automorphism in $\mathrm{Gal}(k/\Q)$. Let $\hat\chi$ be a primitive Hecke character of $k$ of conductor $\mathfrak{f}$ with~$\mathfrak{f}$ an integral ideal. Then $k_{\hchi}/\Q$ is Abelian if and only if $\tau(\mathfrak{f}) = \mathfrak{f}$ and the kernel of~$\hchi$ contains $\mathrm{Cl}_k(\mathfrak{f})^{\tau-1}$.
\end{lemma}

\begin{proof}
Fix $\bar\tau$ an element of $G_\Q$ extending $\tau$. Then $k_{\hchi}^{\bar\tau}$ is an Abelian extension of~$k$ of conductor $\tau(\mathfrak{f})$. Therefore, if $k_{\hchi}/\Q$ is Galois then $\tau(\mathfrak{f}) = \mathfrak{f}$.
Assume that $ \tau(\mathfrak{f}) = \mathfrak{f} $.
Denote by $\tilde{k}(\mathfrak{f})$ the maximal subextension of $k(\mathfrak{f})/k$ that is Abelian over~$\Q$ where $k(\mathfrak{f})$ denotes the ray class field of $\mathfrak{f}$. The extension $k_{\hchi}/\Q$ is Abelian if and only if $k_\hchi \subseteq \tilde{k}(\mathfrak{f})$, that is if and only if the kernel of $\hchi$ contains the subgroup $\tilde{H}$ of $\mathrm{Cl}_k(\mathfrak{f})$ corresponding to $\tilde{k}(\mathfrak{f})$ by class field theory. The Galois group of the extension $k(\mathfrak{f})/\tilde{k}(\mathfrak{f})$ is the commutator subgroup of $\mathrm{Gal}(k(\mathfrak{f})/\Q) = \langle \mathrm{Gal}(k(\mathfrak{f})/k), \bar\tau \rangle$. It is equal to $\langle \bar\tau g \bar\tau^{-1}g^{-1} : g \in \mathrm{Gal}(k(\mathfrak{f})/k) \rangle$ where, by abuse, we still denote by $\bar\tau$ its restriction to $k(\mathfrak{\mathfrak{f}})$. For $g \in  \mathrm{Gal}(k(\mathfrak{f})/k)$, we have $\bar\tau g \bar\tau^{-1} = \Art_\mathfrak{f}(C^\tau)$ where $\Art_\mathfrak{f} : \mathrm{Cl}_k(\mathfrak{f}) \to \mathrm{Gal}(k(\mathfrak{f})/k)$ is the Artin map, and $C \in \mathrm{Cl}_k(\mathfrak{f})$ is such that $g = \Art_\mathfrak{f}(C)$. We conclude that $\tilde{H} = \mathrm{Cl}_k(\mathfrak{f})^{\tau - 1}$ and the lemma is proved. 
\end{proof}

We also use the following reduction.  For $ s $ in $ \Zp $, it is immediate from the construction of $ p $-adic Hecke $ L $-functions that $\Lp p k {\hchi} s $ takes value in $E = \Qp(\hchi)$, the field generated over $\Qp$ by the values of $\hchi$, and that
$\Lp p k {\hchi^\sigma} s = \Lp p k {\hchi} s {}^\sigma$ for  $\sigma$ in $ \Gal(E/\Qp)$. Therefore,  Conjecture~\ref{mainconjecture} is satisfied for $\hchi^\sigma$ if and only if it is satisfied for $\hchi$. The same is true for Conjectures \ref{Grootsconjecture}, \ref{gross conj}, and \ref{muconjecture}. Therefore it is enough to consider primitive Hecke characters up to conjugation over $\Qp$.

Henceforth we consider primitive, non-trivial $p$-adic Hecke characters $\hchi$
over the field $k = \Q(\sqrt{d})$, where $d=1$ (for $ k=\Q $)
or a positive square-free integer,
such that the conductor of $\hchi$ does not have an infinite part, 
$k_\hchi/\Q$ is not Abelian when $d > 1$,
and up to conjugation over $\Qp$. 
For fixed $d$, we consider for different primes $p$ all the characters satisfying the above conditions and whose conductor is of norm less than a bound $B$.
The precise values of $p$ and $B$ when $ k = \Q $ are given in Table~\ref{tab1}.
The values of $p$ and $B$ when $ k $ is a real quadratic field are listed in Table~\ref{tab2}, 
where we give the value of the discriminant of $k$ rather than the value of $d$.
We also give the number (of conjugacy classes) of characters for each discriminant in Table~\ref{tab3}. 

\begin{table}[ht]
\caption{Values of $p$ and $B$ for $k = \Q$.  \label{tab1}}
\vskip-8pt
\begin{tabular}{c|c||c|c}
$p$ & $B$ & $p$ & $B$ \\
\hline
2, 3, 5, 7 & 4999 & 11, 13  & 3999  \\
17, 19   & 3499 & 23, 29   & 2999 \\
31, 37   & 1499 \\
\end{tabular}
\end{table}

\begin{table}[ht]
\caption{Values of $p$ and $B$ for $k$ a real quadratic field.  \label{tab2}}
\vskip-8pt
\begin{tabular}{c|c|c}
disc.  of $k$ & $p$ & $B$ \\
\hline
5, 8  & 2, 3, 5, 7 &  499 \\
12, 13, 17, 21, 24, 28 & 2, 3, 5, 7 & 299 \\
29, 33, 37, 40, 41, 44 & 2, 3, 5, 7 & 299 
\end{tabular}
\end{table}

\begin{table}[h]
\caption{Number of examples.  \label{tab3}}
\vskip-8pt
\begin{tabular}{c|c||c|c||c|c||c|c||c|c}
disc. & \# & disc. & \# & disc. & \# & disc. & \# & disc. & \# \\
\hline
1 & 1\,285\,104 & 12 & 309 & 21 & 308 & 29 & 388 & 40 & 854 \\
 5 & 460 & 13 & 306 & 24 & 336 & 33 & 393 & 41 & 481 \\
8 & 480 & 17 & 385 & 28 & 498 & 37 & 550 & 44 & 469
\end{tabular}
\end{table}

Assume $p$, $d$ and $B$ are fixed. We now explain the construction of all the corresponding characters, up to conjugation over $\Q_p$. 
We use the methods and results of \cite[Chap.~4]{Coh00}. All the computations were performed using PARI/GP \cite{PARI2} and, in fact, most of the algorithms described in \cite[Chap.~4]{Coh00} are implemented in PARI/GP.  The first step is to list all the integral ideals $\mathfrak{f}$ of $k$ with norm at
most $B$ that are conductors. For each such $\mathfrak{f}$, we compute the ray
class group $\mathrm{Cl}_k(\mathfrak{f})$ and list the subgroups $H$ of $\mathrm{Cl}_k(\mathfrak{f})$. We keep only the subgroups $H$ that are primitive and such that $\mathrm{Cl}_k(\mathfrak{f})/H$ is cyclic. Indeed, there exists a Hecke character of $\mathrm{Cl}_k(\mathfrak{f})$ having $H$ as kernel if and only if $\mathrm{Cl}_k(\mathfrak{f})/H$ is cyclic, and this character is primitive if and only if $H$ is primitive. If $k$ is a real quadratic field and $ \mathfrak{f}$ is such that $ \tau(\mathfrak{f}) = \mathfrak{f} $, where $\tau$ is the non-trivial automorphism in $\Gal(k/\Q)$, we discard the subgroups $H$ containing $\mathrm{Cl}_k(\mathfrak{f})^{\tau - 1}$. For each remaining subgroup $H$, we construct a complex primitive Hecke character of conductor $\mathfrak{f}$ whose kernel is $H$ in the following way. We use the notations and results of~\cite[\S.~4.1]{Coh00}. Assume that the group $G$ is given by generators $\mathcal{G}$ and relations $D_G$, where $\mathcal{G} = (g_1, \dots, g_s)$ is a vector 
of elements of $G$ and $D_G$ is an $s \times s$ matrix with integral coefficients. That is, the following sequence is exact
\begin{equation*}\xymatrix{
\Z^s \ar[r]^{D_G} & \Z^s \ar[r]^{\mathcal{G}} & G \ar[r] & 1
}\end{equation*}
where $\Z^s$ is the group of column vectors, the first map sends $V \in \Z^s$ to $D_G V$,
and the second map sends $V \in \Z^s$ to  $\mathcal{G}(V)$ where $\mathcal{G}(V) = g_1^{v_1} \cdots g_s^{v_s}$ for $V = (v_1, \dots, v_s)^t \in \Z^s$. Note that we can always assume that $D_G$ is in Smith Normal Form (SNF), that is $D_G$ is a diagonal matrix $\mathrm{diag}(d_1, \dots, d_s)$ with $d_{i+1} \mid d_i$ for $i = 1, \dots, s-1$. We can also assume that $d_s \ne 1$, i.e., the set $\mathcal{G}$ is a minimal set of generators of $G$. The subgroups of $H$ of $G$ are in $1$-to-$1$ correspondence with the integral matrices $M$ in Hermite Normal Form (HNF) such that $M^{-1}D_G$ has integral coefficients. To such a matrix $M$ corresponds the subgroup $H$ with generators $(\mathcal{G}(M_1),\dots,\mathcal{G}(M_s))$, where $M_1,\dots,M_s$ are the columns of $M$, and relations $M^{-1}D_G$. Using the SNF reduction algorithm, we compute two invertible matrices $U$ and $V$ with integral coefficients and such that $UMV = E$ where $E$ is in SNF. Then the quotient group $G/H$ has generators $(\bar\eta_1, \dots, \bar\eta_s)$ and relations~$E$,
where $\mathcal{G}(U^{-1}) = (\eta_1, \dots, \eta_s)$ and the bar $\bar{\ }$ denotes the reduction modulo $H$. Thus $G/H$ is cyclic if and only if $s = 1$ (that is $G$ itself is cyclic) or $e_2 = 1$ where $E =  \mathrm{diag}(e_1, \dots, e_s)$. We assume now that $G/H$ is cyclic and let $n = e_1$ be its order. Let $\zeta_0$ be a primitive $d_1$-th root of unity. The group $\hat{G}$ of complex characters of $G$ has generators $\mathcal{C} = (\psi_1, \dots, \psi_s)$ and relations $D_{\hat{G}} = \mathrm{diag}(d_1, \dots, d_s)$ with
\begin{equation*}
\psi_i(g_j) = 
\begin{cases}
\zeta_0^{f_i} & \text{ if } i = j, \\
1 & \text{otherwise}
\end{cases}
\end{equation*}
where $f_i = d_1/d_i$. 
Let $F = \mathrm{diag}(f_1, f_2, \dots, f_s)$ and set
\begin{equation*}
C =  (d_1/e_1, 0, \dots, 0) U F^{-1}.
\end{equation*}
and $\hchi = \mathcal{C}(C) \in \hat{G}$. For $g = \mathcal{G}(X) \in G$ with $X \in \Z^s$, we have
\begin{equation*}
\hchi(g) = \zeta_0^{C F X}
\end{equation*}
and thus $\hchi(\eta_1) = \zeta$ and $\hchi(\eta_i) = 1$, for $i = 2, \dots, s$, where $\zeta = \zeta_0^{d_1/n}$ is a primitive $n$-th root of unity. Thus, the kernel of $\hchi$ is $H$ as desired. 

Recall that we want to get exactly one character for each class of conjugation of characters of kernel $H$ under the action of $\Gal(\Qbar_p/\Q_p)$. The construction of the character $\hchi$ depends on the value of the primitive $n$-th root of unity $\zeta$. In the complex case, the different characters obtained for the different choices of $\zeta \in \C^\times$ form a complete  class of conjugation of Hecke characters under the action of $\Gal(\Qbar/\Q)$ since the group $\Gal(\Qbar/\Q)$ acts transitively on the set of (complex) primitive $n$-th roots of unity. 
The situation is different in the $p$-adic case. In this case, let $\Phi_n(X)$ be the $n$-th cyclotomic polynomial and let $F_1(X), \dots, F_t(X)$ be the irreducible factors of the factorization of $\Phi_n(X)$ in $\Q_p[X]$. Then the set of $p$-adic primitive $n$-th roots of unity splits under the action of $\Gal(\Qbar_p/\Q_p)$ into $t$ disjoint subsets; each one being the set of roots of one of the factor $F_i$. Therefore, the set of $p$-adic Hecke characters with kernel equal to $H$ splits under the action of $\Gal(\Qbar_p/\Q_p)$ into $t$ different classes and a representative of each class is obtained by choosing, for $i = 1, \dots, t$, the value of $\zeta$ to be one of the roots of $F_i(T)$ in $\Qbar_p$. Therefore, the computation of a system of representatives of the conjugacy classes of the $p$-adic characters with kernel $H$ over $\Qp$ boils down to the computation of the  irreducible factors $F_1(X), \dots, F_t(X)$ in $\Qp[X]$ of the cyclotomic polynomial $\Phi_n(X)$. Using the fact that the factorization of $\Phi_n$, when $p \nmid n$, is easily computed from the factorisation of $\Phi_n$ modulo $p$ thanks to Hensel's lemma, and that $\Phi_p$ is irreducible over $\Qp$, one can compute these factors recursively by means of the following lemma.

\begin{lemma}
Let $n \geq 1$ be an integer. Assume that $\Phi_n(X) = F_1(X) \cdots F_t(X)$ is the factorization of $\Phi_n$ over $\Qp$. Then
$\Phi_{np}(X)  = G_1(X) \cdots G_t(X)$ is the factorization of $\Phi_{np}(X)$ over $\Qp$ where
\begin{equation*}
 G_i(X) =
 \begin{cases}
 F_i(X^p) & \text{ if } p \mid n, \\
 F_i(X^p)/F_i(X) & \text{ if } p \nmid n. \qed
\end{cases}
\end{equation*}
\end{lemma}

\smallskip

We now turn to the numerical verification of the conjectures. Let $\hchi$ be one of the Hecke characters constructed by the above method and let $\chi$ be the associated $1$-dimensional Artin character. 
We explain the method we use to verify Conjecture~\ref{gross conj}.
Denote by $r_p(\chi)$ the order of vanishing at $s=0$ of the truncated complex $L$-function
\begin{equation*}
  L(k , \sigma \circ\chi\o_p^{-1} , s) \prod_{\mathfrak{p} \in P}  \Eul_\mathfrak{p}(k,\sigma \circ \chi\o_p^{-1},s) 
\,,
\end{equation*}
where $ P $ is the set of primes ideals of  $ k $ lying above~$p$ and $\sigma : \Qpbar \to \mathbb{C}$ is any embedding. By \cite[(2.8)]{Gross81}, $r_p(\chi)$ is equal to the number of prime ideals $\mathfrak{p} \in P$ such that $D_\mathfrak{p}$ is included in the kernel of $\sigma \circ \chi\omega_p^{-1}$ where $D_\mathfrak{p}$ is a decomposition group of $\mathfrak{p}$ in $G_k$; note that the condition does not depend upon the choice of the decomposition group $D_\mathfrak{p}$. Since the kernel of $\sigma \circ \chi\omega_p^{-1}$ is equal to the kernel of $\chi\omega_p^{-1}$, we have $r_p(\chi) = \#R_p(\chi)$ where
\begin{equation*}
R_p(\chi) = \{\mathfrak{p} \in P : \chi(\delta) = \omega_p(\delta),\ \forall \delta \in D_\mathfrak{p} \}. 
\end{equation*}

Let $\hchi$ be the primitive Hecke character of $k$ associated to $\chi$, and denote by $\mathfrak{f}$ its conductor. Let $\hat\omega_p$ be the primitive Hecke character of $k$ associated to $\omega_p$; its conductor is $\mathfrak{wz}$ where $\mathfrak{w}$ is a divisor of $(q)$ and $\mathfrak{z}$ is the set of infinite (real) places of $k$. In our situation $\mathfrak{w}$ is easy to compute thanks to the following result whose proof is left to the reader. 
\begin{lemma}
Let $k = \Q(\sqrt{d})$ be a real quadratic field where $d > 1$ is a square-free integer. For $\mathfrak{p} $ in $P$, we have that $\mathfrak{p}$ divides $\mathfrak{w}$ if and only if one of the following cases applies:
\begin{itemize}
\item $p>3$;
\item $p = 3$ and $3 \nmid d$;
\item $p = 2$ and $d \equiv 2 \pmod{4}$. 
\end{itemize}
\end{lemma}
Let $\mathfrak{m}$ be the lcm of $\mathfrak{f}$ and $\mathfrak{w}$. Thus, by composing with the natural surjections $\mathrm{Cl}_k(\mathfrak{mz}) \twoheadrightarrow \mathrm{Cl}_k(\mathfrak{f})$ and $\mathrm{Cl}_k(\mathfrak{mz})  \twoheadrightarrow \mathrm{Cl}_k(\mathfrak{wz})$ respectively,  we can and will, from now on, consider $\hchi$ and $\hat\omega_p$ as (not necessarily primitive) Hecke characters on $\mathrm{Cl}_k(\mathfrak{mz})$.  Thus, we have  
\begin{equation*}
R_p(\chi) = \{\mathfrak{p} \in P : \hchi(C) = \hat\omega_p(C),\ \forall C \in D_\mathfrak{p}(\mathfrak{mz}) \}
\end{equation*}
where $D_\mathfrak{p}(\mathfrak{mz}) \subset \mathrm{Cl}_k(\mathfrak{mz})$ is the inverse image of $D_\mathfrak{p}$ by the Artin map. Fix $\mathfrak{p} \in P$. Write $\mathfrak{m} = \mathfrak{m}_0 \mathfrak{p}^a$ with $a \geq 0$ an integer and $\mathfrak{p} \nmid \mathfrak{m}_0$. Then $D_\mathfrak{p}(\mathfrak{mz})$ is the inverse image by the canonical surjection $\mathrm{Cl}_k(\mathfrak{mz}) \twoheadrightarrow \mathrm{Cl}_k(\mathfrak{m}_0\mathfrak{z})$ of the subgroup of $\mathrm{Cl}_k(\mathfrak{m}_0\mathfrak{z})$ generated by the class of $\mathfrak{p}$.
This provides us with an efficient way to compute the value of $r_p(\chi)$. 
In all our examples, we found that $r_p(\chi) \leq 1$ and thus Conjecture~\ref{gross conj} follows in these cases by a result of Gross \cite[Proposition~2.3]{Gross81}.  

\begin{theorem}\label{th:proofgoss=0}
Conjecture~\ref{gross conj} is satisfied
for every $1$-dimensional $ \Qpbar $-valued even Artin character $\chi$ of a real quadratic field $k$
for the $ p $, $ k $, and $ \chi $ with conductor of norm at most~$ B$ such that $ k_\chi/\Q $ is non-Abelian, as listed in Table~\ref{tab1}. 
\end{theorem}

\begin{remark} \label{grossremark}
Although we have $r_p(\chi) \leq 1$ in all our examples, there exist examples of $1$-dimensional $ \Qpbar $-valued even Artin characters $\chi$ of real quadratic fields $k$ such that $ k_\chi/\Q $ is non-Abelian and for which $r_p(\chi) = 2$. To find such an example, start with a totally complex quartic field $K$ that is not Galois over $\Q$ and contains a real quadratic field $k$. Let $p$ be a prime number that is totally split in $K/\Q$ and let~$\psi : G_k \to \Qbar_p$ be the only non-trivial Artin character  of $k$ whose kernel contains~$G_K$. Then the Artin character $\chi = \psi \o_p : G_k \to \Qbar_p$ is an example of such a character. For example, for $K = \Q(\alpha)$, with $\alpha$ a root of the polynomial $X^4 - X^3 - 6X^2 - 8X + 64$, which has real quadratic subfield $k = \Q(\sqrt{89})$, and~$ p = 2 $, we obtain as $\chi$ the unique character of $G_k$ of order $2$ and conductor $4\mathfrak{p}$ where~$\mathfrak{p}$ is the prime ideal above $17$ generated by $28+3\sqrt{89}$. 

As we discussed in Remark~\ref{zeroremark}, Spie\ss\ \cite{Spiess14} has proved that the order of vanishing of $\Lp p k {\chi} s $ at $s=0$ is at least $r_p(\chi)$. Thus,
in our example $u-1$ is a root of order at least $r_2(\chi)=2$ of $G_\chi(T)$, where~$u$ is the unit such that  \eqref{iwasawaseries} applies. We computed $G_\chi(T)$ and checked numerically that $G^{(3)}_\chi(T)$ does not vanish at~$T = u-1$. Therefore, the order of vanishing of  $\Lp 2 k {\chi} s $ at $s = 0$ is indeed equal to $r_2(\chi)$ and Conjecture~\ref{gross conj} is satisfied for this example. 
\end{remark}

We now explain the methods used to numerically investigate Conjecture~\ref{mainconjecture}.
We start with the following result.
\begin{proposition}\label{zeroesofL} 
 Let $w = v_p(u-1)$ . Then, the map $s \mapsto u^{1-s} - 1$ yields a bijection between the zeroes of $\Lp p k {\chi} s $ in~$\Z_p$ and the roots of $G_\chi(T)$ in $p^w\Z_p$. 
\end{proposition}

\begin{proof}
Let $\exp_p$ and $\log_p$ denote respectively the $p$-adic exponential and logarithmic functions. Observe that, since $\log_p(u) \in p^w\Z_p$, we have $u^{1-s} = \exp_p((1-s)\log_p(u))$ for any $s \in  \Z_p$. By the properties of the exponential and logarithmic functions, the maps $s \mapsto u^{1-s} - 1$ and $t \mapsto 1 - \log_p(1+t)/\log_p(u)$ are inverse bijections between $\Z_p$ and $p^w\Z_p$. Furthermore, by \eqref{iwasawaseries}, $a$ is a zero of $\Lp p k {\chi} s $ if and only if $u^{1-a}-1$ is a root of $G_\chi(T)$. The result is proved. 
\end{proof}

Using the methods of \cite{Rob15} and the fact that $\Lp p k {\chi} s = \Lp p k {\hchi} s $, we compute a polynomial $\tilde{G}_\chi(T)$ with coefficients in $\Zp[\chi]$ such that
\begin{equation*}
\tilde{G}_\chi(T) \equiv G_\chi(T) \pmod{p^{M_0}, T^L}
\end{equation*}
for some positive integers $M_0$ and $L$. In fact, we start with $L$ small, say $L = 3$, and $M_0$ quite big, say such that $p^{M_0} \approx 10^{10}$ with the extra condition that $M_0 \geq w(L+1) + 3$ if $r_p(\chi) = 1$. 
If necessary, we increase the value of $ L $ until we find a coefficient in $ \tilde{G}_\chi(T) $ with $p$-adic valuation as small as possible according to the ``$\mu=0$'' conjecture, i.e., with valuation 0 when $ p $ is odd, and $ [k:\Q] $ when $ p=2 $. When we increase the value of $L$, we also increase that value of $M_0$ if needed in the case where $r_p(\chi) = 1$. We find that it is always possible to find such a value of $L$ and thus we proved in passing the following result on the ``$\mu=0$'' conjecture.
(We recall that this conjecture was proved when $ k = \Q $; see \cite{Fer-Wash}.)

\begin{theorem}\label{th:proofmu=0}
Conjecture~\ref{muconjecture} is satisfied
for every $1$-dimensional $ \Qpbar $-valued even Artin character $\chi$ of a real quadratic field $k$
for the $ p $, $ k $, and $ \chi $ with conductor of norm at most~$ B$ such that $ k_\chi/\Q $ is non-Abelian, as listed in Table~\ref{tab1}. 
\end{theorem}

It follows from Proposition~\ref{zeroesofL}  that the non-zero zeroes of the $p$-adic $L$-function of $\chi$ in $\Zp$ correspond to the roots of $G_\chi(T)$ in $p^w\Zp$ distinct from $u-1$. Set
\begin{equation*}
G^\sharp_\chi(T) = c_p \frac{G_\chi(T)}{(T-(u-1))^{r_p(\chi)}}
\end{equation*}
with $c_p = 1$ if $p$ is odd and $c_p = 2^{-[k:\Q]}$ if $p = 2$. Since Conjecture \ref{gross conj} is satisfied for all the characters $\chi$ that we consider (see Theorem~\ref{th:proofgoss=0}), $G^\sharp_\chi(T)$ does not vanish at $u-1$ and the roots of $G^\sharp_\chi(T)$ are exactly the roots of $G_\chi(T)$ distinct from $u-1$. Furthermore, since
\begin{equation*}
\frac{1}{(T-(u-1))} = \sum_{n \geq 0} (u-1)^{-(n+1)} T^n 
\end{equation*}
we can deduce from our approximation $\tilde{G}_\chi(T)$ of $G_\chi(T)$, an approximation $\tilde{G}^\sharp_\chi(T)$ of $G^\sharp_\chi(T)$ such that $\tilde{G}^\sharp_\chi(T) \equiv G^\sharp_\chi(T) \pmod{p^M, T^L}$ with $M = M_0 - [k:\Q]v_2(p)$ if $r_p(\chi) = 0$ and $M = M_0 - wL - [k:\Q]v_2(p)$ if $r_p(\chi) = 1$. In any case, we always $M > w$, though it is much larger in general.\footnote{The coefficients of $\tilde{G}^\sharp_\chi(T)$ are actually obtained with decreasing $p$-adic precision; although we use this fact in our computations, we will not take it into account here to simplify the exposition.}

We compute approximations modulo~$p^M$ of the roots of $\tilde{G}^\sharp_\chi$ in $p^w\Z$. In many cases, we find that it has no such root and thus $G_\chi(T)$ has no root in $p^w\Z$ distinct from $u-1$ and Conjecture~\ref{mainconjecture} is satisfied for the character $\chi$. Otherwise, we use the following result which is a direct application of Hensel lifting. 

\begin{lemma}
Let $r \in p^w\Z$ be such that $\tilde{G}^\sharp_\chi(r) \equiv 0 \pmod{p^M}$ and let $m$ be the $p$-adic valuation of $(\tilde{G}^\sharp_\chi)'(r)$. Assume that $m < M/2$. Then there exists a unique root $t$ of $G^\sharp_\chi(T)$ in $p\Zp$ such that $t \equiv r \pmod{p^{M - m}}$. \qed
\end{lemma}

We compute the elements $r_i \in p^w\Z$ such that $\tilde{G}^\sharp_\chi(r_i) \equiv 0 \pmod{p^M}$. For each element $r_i$, let $m_i$ be the $p$-adic valuation of $(\tilde{G}^\sharp_\chi)'(r_i)$. If $m_i < M/2$ for all $i$'s, then we store the $r_i$'s as the approximations to the precision $p^{M-m_i}$ of the roots of $G^\sharp_\chi(T)$ in $p^w\Zp$. If $m_i \geq M/2$ for some $i$, then we increase the precision $M$ (and $L$ if necessary) and start over. Eventually, we get to a large enough precision $M$ such that all the elements $r_i \in p^w\Z_p$ such that $\tilde{G}^\sharp_\chi(r_i) \equiv 0 \pmod{p^M}$ satisfy $m_i < M/2$ and thus we deduce approximations of all the roots of $G_\chi(T)$ in $p^w\Z_p$ distinct from $u - 1$. Note that, for this method to work and not end up in an infinite loop, it is necessary that the roots of $G_\chi(T)$ in $p\Zp$, distinct from $u-1$, are simple. We find that it is indeed the case in all the examples that we consider. (In fact, we have the stronger statement that Conjecture~\ref{Grootsconjecture} holds for all the characters that we consider, see Theorem~\ref{th:simple}.)

In the cases where we find that $G^\sharp_\chi(T)$ admits some root in $p^w\Z_p$, and thus $G_\chi(T)$ admits some root in $p^w\Z_p$ distinct from $u-1$, we cannot decide this way if Conjecture~\ref{mainconjecture} is satisfied or not. Indeed, it is impossible computationally to distinguish between the approximation of a $p$-adic integer and a rational integer. Thus,  we need to convince ourselves that these roots are actually approximations of elements of $\Z_p$ that do not lie in $\Z$. For that, we use the following criterion. For $x \in \mathbb{Z}_p$ and $E \geq 1$, denote by $(x \!\mod p^E)$ the unique non-negative integer less than $p^E$ such that $x - (x \mod \!p^E) \in p^E\Z_p$. Assume that $a \in \Z_p$, $a \not= 0$, is a zero of $\Lp p k {\chi} s $. Then $a$ does not lie in $\Z$
if and only if the $p$-adic expansion $a = \sum_{i \geq 0} a_i p^i$, with $a_i \in \{0, \dots, p-1\}$, is infinite (since, if $a \in \Z$ then $a \geq 0$ as
observed in Section~\ref{nozeroessection}).
In the computations, once we have computed an approximation $(a \!\mod p^E)$ of the zero $a \in \Z_p$ by the above method, we compute in the same way a new approximation $(a \!\mod p^{E'})$ with $E' = E+3$. Assuming that $a \not\in \mathbb{Z}$ and that the coefficients $a_i$ are randomly distributed, we expect $(a \!\mod p^E) \not= (a \!\mod p^{E'})$ with probability $1 - 1/p^3 \geq 87\%$ (and in fact $\geq 96\%$ for $p \geq 3$). If this is not the case, we increase the value of $E'$ and test the condition again. We have verified in all our examples that indeed $(a \!\mod p^E) \not= (a \!\mod p^{E'})$ for $E' = E+3$ in most cases and, in some few cases, for $E' = E + 6$. Of course, as mentioned above, this does not prove that $a$ is not a rational integer as the same thing could happen if $a$ was an integer greater than $p^E$. It does show however that any non-zero zero of $\Lp p k {\chi} s $ is either an element of $\Z_p \setminus \Z$ or a huge positive integer.

\begin{theorem}\label{th:mainconj}
Let $ k = \Q $ or a real quadratic field, and $ \chi $
a $ 1 $-dimensional~$ \Qpbar $-valued even Artin character of
$ k $ such that $ k_\chi/\Q $ is non-Abelian if $ k \ne \Q $.
\begin{enumerate}
\item
If $ k = \Q $, and $ (p,B) $ occurs in Table~\ref{tab1},
then the function $ \Lp p k {\chi} s $ has no root $ m = 1, 2, \dots, 10^8 $
if the conductor of $ \chi $ is at most $ B $.

\item
If $ k $ is a real quadratic field such that it and $ (p,B) $ occur in Table~\ref{tab2},
then the function $ \Lp p k {\chi} s $ has no root $ m = 1, 2, \dots, 10^8 $
if the norm of the conductor of $ \chi $ is at most $ B $.
\end{enumerate}

Furthermore, for $1\,283\,351$ of the $1\,291\,321$ conjugacy
classes over $\Q_p$ of the $\chi$ considered here,
the function $\Lp p k {\chi} s $ has no non-zero zero in
$\Z_p$ and thus satisfies Conjecture~\ref{mainconjecture}.
\end{theorem}

\begin{remark}
The precise number of classes considered for which the $L$-function does not have a non-zero zero in $\Z_p$ for each field discriminant is listed in Table~\ref{tab4}. 
\end{remark}

\begin{table}[ht] 
\caption{The number of classes with no non-zero zero in $\Z_p$. \label{tab4}}
\vskip-8pt
\begin{tabular}{c|c||c|c||c|c||c|c||c|c}
disc. & \# & disc. & \# & disc. & \# & disc. & \# & disc. & \# \\
\hline
1 &  1\,277\,552 & 12 & 276 & 21 & 288 & 29 & 367 & 40 & 803 \\
5 & 438 & 13 & 292 & 24 & 314 & 33 & 357 & 41 & 444 \\
8 & 440 & 17 & 371 & 28 & 460 & 37 & 518 & 44 & 431
\end{tabular}
\end{table}

\begin{remark}\label{rk:schneider}
As noted in Section~\ref{consequences}, part of a conjecture by Schneider (see \cite[p.~192]{schn79}) holds
if, for $ p $ odd, $ k $ totally real and $ m $ even, one has $ \Lp p k {\o_p^m} 1-m \ne 0 $
for $ m \ne 0 $ and $ \zeta_p(k,s) $ having a pole of order~1 at $ s = 1 $ for $ m = 0 $.

By Frobenius reciprocity we have $ \Ind {G_k} {G_\Q} \o_p^m = \o_p^m + \o_p^m \chi_k $
with $ \chi_k $ the non-trivial character of $ G_\Q $ that is
trivial on $ G_k $, so $ \Lp p k {\o_p^m} s = \Lp p {\Q} {\o_p^m} s  \Lp p {\Q} {\o_p^m \chi_k} s $,
and this is $ \Lp p {\Q} {\chi_k} s  \zeta_p (k,s) $ precisely
when one of $ \o_p^m $ and $ \o_p^m \chi_k $ equals~$ {\bf 1}_\Q $.
So from Lemma~\ref{vlemma} we see that $ \Lp p k \o_p^m s $ has the desired behaviour
for $ s = m $ if
and only if $ \Lp p {\Q} {\o_p^m} 1-m \ne 0  $ if $ \o_p^m \ne {\bf 1}_\Q $
as well as  $ \Lp p {\Q} {\o_p^m\chi_k} 1-m \ne 0  $ if $ \o_p^m\chi_k \ne {\bf 1}_\Q $.

By~\eqref{interpol} we have that $ \Lp p {\Q} {\chi} m \ne 0 $ for any
even 1-dimensional character $ \chi $ of $ G_\Q $ when $ m < 0 $,
and it holds for many $ p $, $ \chi $ and $ m > 0 $ by Theorem~\ref{th:mainconj}(1),
and for many $ p $, $ \chi $ and all $ m $ by the last statement
of that theorem. Therefore, our calculations prove many instances of Schneider's conjecture.
\end{remark}

Finally, we turn to proving Conjecture~\ref{Grootsconjecture}
in many cases by means of computer calculations. First, we prove that the conjecture applies to all the Artin characters that we consider.

\begin{proposition}
Let $\chi$ be a $1$-dimensional even Artin character of $k$ with $k = \Q$ or a real quadratic field,
such that $k_\hchi/\Q$ is non-Abelian if $k \not= \Q$. Then $\Ind {G_k} {G_\Q} (\chi)$ is irreducible.
\end{proposition}

\begin{proof}
The result is clear if $k = \Q$. Assume now that $k$ is a real quadratic field. Let $\chi' = \Ind {G_k} {G_\Q} (\chi)$ and assume that $\chi'$ is reducible. Since $\chi'$ is of dimension $2$, it is the sum of two $1$-dimensional characters, say $\chi' = \nu_1 + \nu_2$,
with $ \nu_1 \nu_2^{-1} = \chi_k $, the unique character of $G_\Q$ with $G_k$, and the $ \nu_i $ restricting to $ \chi $ on $ G_k $.
Therefore $  \ker(\chi) = \ker(\nu_1) \cap \ker(\chi_k) = \ker(\nu_1) \cap \ker(\nu_2) = \ker(\chi') $.
But the $G_\Q/\ker(\nu_i)$ are Abelian (in fact cyclic), so $G_\Q/\ker(\chi') = G_\Q/\ker(\chi)$ is also Abelian, a contradiction. Hence $\Ind {G_k} {G_\Q} (\chi)$ is irreducible.
\end{proof}

The statement of Conjecture~\ref{Grootsconjecture} is true for the character $\chi$ if and only if the series $G^\sharp_\chi(T)$ has only simple roots of positive valuation. Fix a valuation ring $\O$ in a finite extension of $\Qp$ such that $G^\sharp_\chi(T)$ lies in $\O[[T]]$. In our computations, we always take $\O$ to be the smallest valuation ring satisfying this condition. Let $\pi$ be a uniformizing parameter of $\O$ and let $\kappa = \O/\pi\O$ be its residue field. The series $G^\sharp_\chi(T)$ has only simple roots of positive valuation if at least one of the two following conditions holds:
\begin{enumerate}
\item[(G1)] $\lambda(G^\sharp_\chi) = 0$ or $1$;
\item[(G2)] $v_p(G^\sharp_\chi(0)) = v_p(\pi)$.
\end{enumerate}
Condition (G1) is clear since in this case $G^\sharp_\chi$ has at most one root of positive valuation. For condition (G2), the distinguished polynomial $D^\sharp_\chi(T)$ of $G^\sharp_\chi(T)$ lies in $\O[T]$ and, since $v_p(D^\sharp_\chi(0)) = v_p(G^\sharp_\chi(0)) = v_p(\pi)$, it is an Eisenstein polynomial at $\pi$; in particular, it is irreducible and thus has only simple roots. Both conditions are easy to check on the computed approximation $\tilde{G}^\sharp_\chi$ of $G^\sharp_\chi$. 

These two conditions are not enough in general to prove that the roots of positive valuation of $G^\sharp_\chi(T)$ are simple. For the remaining cases, we use conditions on the Newton polygon of $G^\sharp_\chi(T)$ and possibly that of its derivative $(G^\sharp_\chi)'(T)$. Indeed, since Conjecture~\ref{muconjecture} is satisfied for all the characters $\chi$ that we consider (see Theorem~\ref{th:proofmu=0}), the power series $G^\sharp_\chi(T)$ has integral coefficients and at least one coefficient with zero valuation, so that the methods of Section~\ref{NPsection} apply. Information on these Newton polygons is deduced from the Newton polygons of the computed approximations $\tilde{G}^\sharp_\chi(T)$ of  $G^\sharp_\chi(T)$, and its derivative, using Lemma~\ref{lem:incnp}.

The roots of positive valuation of the power series $G^\sharp_\chi(T) $ are simple if, for every segment $\sigma \in \mathcal{NP}(G^\sharp_\chi)$, at least one of the following conditions holds,
where the polynomial $P_\sigma$ in (S1) and (S3) is defined in Lemma~\ref{linklemma}:
\begin{enumerate}
\item[(S1)] $P_\sigma(T) = T^{l(\sigma)} + u$ with $u \not= 0$ and $p \nmid l(\sigma)$, where $ l(\sigma) $ is the length of $ \sigma $;
\item[(S2)] the integers $l(\sigma)$ and $e h(\sigma) $ are relatively prime where $l(\sigma) $ and $h(\sigma) $ are the length and the height of $\sigma$ and $e \geq 1$ is such that $v_p(\pi) = 1/e$;
\item[(S3)] $p = 2$ and $P_\sigma(T) = T^2 + u_1T + u_0$ with $u_0, u_1 \not= 0$;
\item[(S4)] there is no segment in $\mathcal{NP}((G^\sharp_\chi)')$ of slope equal to $s(\sigma)$;
\item[(S5)] there exist two polynomials $A, B \in \O[T]$, with $B \not = 0$, such that there is no segment in $\mathcal{NP}(AG^\sharp_\chi+B(G^\sharp_\chi)')$ of slope equal to $s(\sigma)$.  
\end{enumerate}
Observe that (S4) is a special case of (S5), with $A = 0$ and~$B = 1$. However,
we specify both conditions in order to emphasise the fact that we do test (S4) before testing~(S5). The condition (S5) is used as a last resort when all the other conditions fail.

\begin{figure}[htp]
\begin{tikzpicture}[scale=0.38]
\draw (0,8) -- (0,0) -- (32,0);
\foreach \Point in {(0,7), (1,3), (2,1), (3,2), (4,1), (5,3), (6,3), (7,2), (8,3), (9,2), (10,2), (11,5), (12,2), (13,2), (14,3), (15,2), (16,1), (17,1), (18,2), (19,3), (20,4), (21,2), (22,5), (23,1), (24,1), (25,1), (26,2), (27,1), (28,1), (29,1), (30,1), (31,0)}{ \node at \Point {\textbullet};}; 
\draw (0,7) -- (1,3) -- (2,1) -- (31,0); 
\end{tikzpicture}
\caption{An example of Newton polygon ($k = \Q$, $p=2$, conductor $508$, order $2$)}\label{fig:3}
\end{figure}
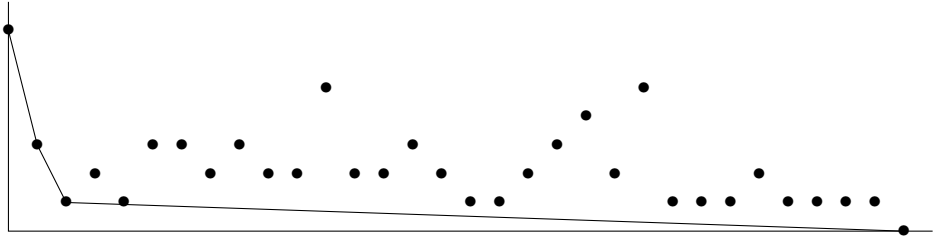

We justify the different conditions.

If condition (S1) is satisfied, then $P_\sigma(T)$ is square-free and thus all the roots of~$G^\sharp_\chi(T)$ of valuation $-s(\sigma)$ are distinct.

For condition~(S2), we shall prove that the distinguished polynomial $D_s(T)$ corresponding to~$ \sigma $ is irreducible in $\O[T]$,
where $ s = h(\sigma) / l (\sigma) $, so that $G^\sharp_\chi(T)$
can have no multiple roots of valuation $-s$.
Normalizing $ v_p(\pi) = 1/e $ to $ 1 $, we see that any such
root has normalized valuation $ -se = eh(\sigma) / l(\sigma) $,
so that by our assumption it generates a field extension over the fraction field
of $ \O $ of degree divisible by $ l(\sigma) $. As it is a root
of $ D_s(T) $, which has degree $ l(\sigma) $, it follows that
$ D_s(T) $ is irreducible.

For condition (S3), one verifies readily that any polynomial over a finite field of characteristic $2$ of the form $T^2 + u_1T + u_0$, with $u_0, u_1 \not= 0$,  is square-free.

For conditions (S4) and (S5), observe that any multiple root of $G^\sharp_\chi(T)$ of valuation $-s(\sigma)$ is also a root of $(G^\sharp_\chi)'(T)$ and thus $\mathcal{NP}(AG^\sharp_\chi+B(G^\sharp_\chi)')$ contains a segment of slope $s(\sigma)$.

The polynomials we use for (S5) are constructed in the following way. Fix two integers $M_1$ and $L_1$ with $0 < M_1 \leq M$ and $0 < L_1 \leq L$. Let $G_1$ be a polynomial in $\O[T]$ such that $G_1 \equiv G^\sharp_\chi \pmod{p^{M_1}, T^{L_1}}$. Using, for example, the sub-resultant algorithm, we compute two polynomials $A$ and $B$ such that
\begin{equation*}
A G_1 + B G_1' = \delta
\end{equation*}
for some $\delta \in \O$ with minimal valuation. Then, we check to see if $\mathcal{NP}(AG^\sharp_\chi+B(G^\sharp_\chi)')$ contains a segment of slope equal to $s(\sigma)$. As noted above, since we only know an approximation of $G^\sharp_\chi(T)$, we cannot always deduce the full Newton polygon $\mathcal{NP}(AG^\sharp_\chi+B(G^\sharp_\chi)')$. Still, using Lemma~\ref{lem:incnp} we can get enough information on the slopes of $\mathcal{NP}(AG^\sharp_\chi+B(G^\sharp_\chi)')$ to conclude in many cases. In the remaining cases where we cannot find by this method two polynomials $A$ and $B$ for which (S5) applies, we resort to choosing random polynomials $A$ and $B$. In a few examples, we could not find after a reasonable amount of time two suitable polynomials $A$ and $B$. Quite often, these were examples for which the $\lambda$-invariant $\lambda(G_\chi)$ is slightly less than the precision $L$ in $T$. For these examples, we recomputed approximations of $G_\chi(T)$ and $G^\sharp_\chi(T) $  to a larger precision~$L$ in $T$ and tested again. In all cases, eventually, we could conclude that all the roots of positive valuation of $G^\sharp_\chi(T) $ in $\Qpbar$ are simple and thus obtained the proof of the following theorem.

\begin{theorem} \label{th:simple}
Let $ k = \Q $ or a real quadratic field, and $ \chi $
a $ 1 $-dimensional~$ \Qpbar $-valued even Artin character of
$ k $ such that $ k_\chi/\Q $ is non-Abelian if $ k \ne \Q $.
\begin{enumerate}
\item
If $ k = \Q $, and $ (p,B) $ occurs in Table~\ref{tab1},
then the character $\chi$ satisfies Conjecture~\ref{Grootsconjecture} if the conductor of $\chi$ is at most $B$. 

\item
If $ k $ is a real quadratic field such that it and $ (p,B) $ occur in Table~\ref{tab2},
then  the character $\chi$ satisfies Conjecture~\ref{Grootsconjecture} if the norm of the conductor of $ \chi $ is at most $ B $.
\end{enumerate}
\end{theorem}

%% file: lambdamodel.tex
\section{Numerical study of $\lambda$-invariants of \padic \Lfunctions} \label{lambdamodel}

From the large set of Iwasawa power series that we computed 
as described in the previous section, 
we tried to understand the behaviour of the resulting coefficients
in the various $ \Z_{p,\chi} $. Our first thought was that they should
behave essentially randomly, but this hope was quickly dashed by
the appearance of distinguished polynomials of, for such random
behaviour, extraordinarly large degree (see, e.g., Example~\ref{largeexample}).
But those exceptions all had the conspicuous
property that the order of the character $ \chi $ was divisible by~$ p $.
In order to understand the behaviour in this case, we below derive
Corollary~\ref{lambdappowercase} from  a special version 
of the Riemann-Hurwitz genus formula of Kida due to Sinnott.
We quote this essential tool for the study of~$\lambda$-invariants of characters of order divisible by~$p$
as Theorem~\ref{sinnottthm} below.

After imposing the condition that $ \chi $ is of order not divisible
by~$ p $, the behaviour appears to be more random. In order to have enough
data available and simplify the setup, we consider only the constant term in the Iwasawa
power series and always use $ k = \Q $. More precisely, we consider the case when the 
constant term is in the maximal ideal which is equivalent to the fact that the 
$\lambda$-invariant of $\chi$ is positive. 
Under our hypothesis that the behaviour of this constant coefficient should be random and since the extension $\Q_p(\chi)/\Q_p$ is unramified, we expected therefore the
$\lambda$-invariant of $\chi$ to be positive with probability $p^{-[\Q_p(\chi):\Q_p]}$. 

In this case, we found that in the numerical data the 
probability that $\lambda(\chi)$ is positive was in general significantly higher
than predicted by our model.
However, after investigating, it turned out that, for an odd prime~$ p $, 
if we take the character of the form~$ \o_p^i \psi $ with $ \psi $
of conductor and order prime to~$ p $, and~$ i $ and $ \psi $
both even, then the numerical data appears to be more in accordance with our model. 
This leads us to making Conjecture~\ref{randomconjecture}, which is
 the same as Conjecture~\ref{Qconjecture} formulated in the
introduction.
We corroborate this by means of the graphs presented at the end of this section. 
In the case $p = 2$, or when the character $\chi$ is not of the form mentioned above, we were not able to find a suitable model for when the $\lambda$-invariant of $\chi$ is positive. 

We now investigate the case when the order of $ \chi $
is divisible by $ p $, returning to the `randomness' of the constant
term of the Iwasawa power series after Example~\ref{lastexample}.

The following result is due to Sinnott~\cite[Theorem~2.1]{Sinnott}. 

\begin{theorem} \label{sinnottthm}
Let $k$ be a totally real number field and let $ p $ be any prime number. Let $\theta$ and $\psi$ be two even $1$-dimensional characters of
$\Gal(k^{ab}/k)$ satisfying the ``$\mu = 0$'' conjecture $($see Conjecture~\ref{muconjecture}$)$. Assume that the order of $\psi$ is a power of $p$ and the order of $\theta$ is prime to $p$. Then
\begin{equation*}
\lambda(\theta\psi) = \lambda(\theta) + N(\theta,\psi)
\,,
\end{equation*}
where $N(\theta,\psi)$ is the number of prime ideals of the ring of integers $ \O_{k_\infty} $ of $k_\infty$ that divide the conductor of $\psi$ but not $p$, and are totally split in the compositum of~$k_{\theta\omega_p^{p-2}}$ and $k_\infty$. \qed
\end{theorem}

As a direct application of the result of Sinnott, we give a formula for the $\lambda$-invariant of characters of $p$-power order. 

\begin{corollary} \label{lambdappowercase}
Let $k$ be a totally real number field and let $ p $ be any prime number.
Let $e \geq 0$ be the integer such that $ k \cap \, \Q_\infty = \Q_e$
where $\Q_n$ denotes the $n$-th layer in the $\Z_p$-extension $\Q_\infty/\Q$. For an ideal
$\mathfrak{m}$ of the ring of integer $\O_k$ of $k$, prime to $p$, define 
\begin{equation*}
a_\mathfrak{m} = v_p\Big(\frac{\langle\mathcal{N}(\mathfrak{m})\rangle - 1}{q}\Big)
\end{equation*}
where $\mathcal{N}(\mathfrak{m}) = \mathrm{card}(\O_k/\mathfrak{m})$ is the norm of the ideal $\mathfrak{m}$ and $x \mapsto \langle x \rangle$ is the projection of $\mathbb{Z}_p^\times$ onto $1+q\mathbb{Z}_p$.
Then for any even $1$-dimensional Artin character $\chi : G_k \to \Qpbar $ of conductor $\mathfrak{f}$ and $p$-power order, under the assumption that $\chi$ and the trivial character ${\bf 1}_k$ of $k$ satisfy the ``$\mu = 0$'' conjecture $($see Conjecture~\ref{muconjecture}$)$, we have
\begin{equation*}
\lambda(\chi) = \lambda({\bf 1}_k) + \sum_{\substack{\mathfrak{q} \mid \mathfrak{f} \\ \mathcal{N}(\mathfrak{q}) \in  1+p\Z}} p^{a_\mathfrak{q} - e}
\end{equation*}
where the sum is over the prime ideals $\mathfrak{q}$ of $\O_k$ dividing $\mathfrak{f}$
and $ \mathcal{N}(\mathfrak{q}) $ in $ 1+p\Z $.
\end{corollary}

\begin{proof}
We apply Theorem~\ref{sinnottthm} with $\psi = \chi$ and $\theta = {\bf 1}_k$. Denote by $L$ the compositum of $k_{\omega_p^{p-2}}$ and $k_\infty$ and by $\mathcal{A}$ the set of prime ideals of $\O_{k_\infty}$ that divide the conductor~$\mathfrak{f}$ of $\chi$, are not above $p$, and are totally split in $L$. Therefore $N({\bf 1}_k, \psi) = \#\mathcal{A}$ by the theorem and it remains to compute~$ \#\mathcal{A} $.

Let $\mathfrak{Q}$ be a prime ideal of $\O_{k_\infty}$ and let $\mathfrak{q}$ be the prime ideal of $\O_k$ below $\mathfrak{Q}$. Observe that $\mathfrak{Q}$  divides $\mathfrak{f}$ and is not above $p$ if and only if the same is true for $\mathfrak{q}$.

Also note that  $\mathfrak{Q}$ is totally split in $L/k_\infty$ if
and only if the same holds for all primes of $ k_\infty $ above $ \mathfrak{q} $
as they are conjugate, so that $ \mathfrak{Q} $ is totally split
in $ L/k_\infty $ if and only if the number of primes of $ L $
above $ \mathfrak{q} $ is divisible by $ [L : k_\infty] $ as
this degree does not contain a factor~$ p $.
Using the tower $ L/k(\mu_p)/k $ instead, one sees similarly that this number
of primes is divisible by $ [L : k_\infty] = [k(\mu_p) : k] $
if and only if $ \mathfrak{q} $ splits completely in $ k(\mu_p) $.
So we have to count the number of $ \mathfrak{Q} $ such that
the corresponding $ \mathfrak{q} $ divides $ \mathfrak{f} $, is
not above~$ p $, and splits completely in $ k(\mu_p) $.

A prime $ \mathfrak{q} $ of $ k $ that is not above~$ p $ splits
completely in $ k(\mu_p) $ precisely when $ \mathcal{N}(\mathfrak{q}) \equiv 1 $
modulo~$ p $. In order to see this, note that $ \mu_p $
injects into $ \O_k/\mathfrak{q} $ as $ \mathfrak{q} $ is not
above~$ p $. So if~$ \mathfrak{q} $ splits
completely in $ k(\mu_p) $ then $ \mathcal{N}(\mathfrak{q}) \equiv 1 $
modulo~$ p $.
Conversely, assume $ \mathcal{N}(\mathfrak{q}) \equiv 1 $ modulo~$ p $,
and let $ \zeta_p $ be a primitive $ p $th root of unity.
Then~$ \mathfrak{q} $ splits into $ [k(\mu_p) : k] $ different
ideals in $ \O_k[\zeta_p] \subseteq \O_{k(\mu_p)}$ already because
the minimal polynomial of~$ \zeta_p $ over~$ k $
factors into different linear factors if we reduce its coefficients
from $ \O_k $ to~$ \O_k / \mathfrak{q} $,
hence $ \mathfrak{q} $ splits completely in $ \O_{k(\mu_p)} $.

Because a prime $ \mathfrak{q} $ with $ \mathcal{N}(\mathfrak{q}) \equiv 1 $
modulo~$ p $ is not above~$ p $, we now know that
\begin{equation*}
\mathcal{A} = \bigcup_{\substack{\mathfrak{q} \mid \mathfrak{f} \\ \mathcal{N}(\mathfrak{q}) \in 1+p\Z}} \mathcal{A}(\mathfrak{q}) 
\end{equation*}
where $\mathcal{A}(\mathfrak{q})$ is the set of prime ideals of $k_\infty$ above $\mathfrak{q}$.
Let $n \geq 0$.
Because $ \mathfrak{q} $ does not ramify in $ k_n/k $, the number
of prime ideals of $ k_n $ above $ \mathfrak{q} $ is the index
of the decomposition group of such a prime in $ \Gal(k_n/k) $.
This decomposition group maps isomorphically to the Galois group
of the residue field extension, and the action on~$ k_n(\mu_q) $
of a generator corresponding to the Frobenius is determined
by~$ \mathcal{N}(\mathfrak{q}) $ in~$ \mu_{\phi(q)} \cdot (1+q\Zp) $,
its action on $ k_n $ by $ \langle \mathcal{N}(\mathfrak{q}) \rangle $
in~$ 1+q\Zp $.
(Note that~$ \mathcal{N}(\mathfrak{q}) = \langle \mathcal{N}(\mathfrak{q}) \rangle $
for~$ p \ne 2 $ but that for $ p = 2 $ they may not be the same.)

Since $\Gal(k_n/k) \simeq (1+qp^e\Z)/(1+qp^{e+n}\Z)$, with
the image of our generator corresponding to the class of $\langle \mathcal{N}(\mathfrak{q})\rangle$,
the decomposition group has order~$p^{e+n-a_\mathfrak{q}}$ for~$n \geq a_\mathfrak{q}-e$,
so there are $p^n/p^{e+n-a_\mathfrak{q}} = p^{a_\mathfrak{q}-e}$ prime ideals above $\mathfrak{q}$ in $k_n$ for~$n$ large enough.
Thus $\# \mathcal{A}(\mathfrak{q}) = p^{a_\mathfrak{q}-e}$ and the result is proved. 
\end{proof}

\begin{remark}\label{compl1}
When $k = \Q$, we have $\lambda({\bf 1}_\Q) = -1$ since $\frac{1}{2}G_{{\bf 1}_\Q}(T)$ is a unit of $\Z_p[[T]]$ by \cite[Lemma~7.12]{wash} and $H_{{\bf 1}_\Q} (T)= T$. When $k$ is a real quadratic field, then the $p$-adic zeta function of $k$ is the product of the $p$-adic zeta function of $\Q$ and of the $L$-function $ \Lp p {\Q} {\chi_k} s $ where $\chi_k$ is the non-trivial character of $ G_\Q $ that is trivial on $ G_k $. 
Therefore $\lambda({\bf 1}_\Q+\chi_k) = \lambda(\chi_k) - 1$ and one computes 
 $ \l({\bf 1}_k) $ from this by comparing generators of $ \Gal(k_\infty/k) $
 and $ \Gal(\Q_\infty/\Q) $ as around~\eqref{twoTs}, so that
 $ \l({\bf 1}_k) $ equals $ \lambda(\chi_k) - 1 $ if $e = 0$ or $ \frac12 (\lambda(\chi_k) - 1 ) $ if $e = 1$ 
 (this is the case only if $k = \Q(\sqrt{2})$). 
 Note that we can have  $\lambda({\bf 1}_k) \geq 0$ (see, e.g., Example~\ref{exl1=0}). 
\end{remark}

The following generalises Lemma~\ref{vlemma}(1).

\begin{proposition} \label{vprop}
Let $ p $ be a prime number.
If $ \chi : G_\Q \to \Qpbar $ is an even 1-dimensional Artin character 
of $ p $-power order with conductor $ p^l $ for $ l \ge 1 $, then
$ \l(\chi) = -1 $ and
$ \frac12 G_{\chi}(T) $ is in $ \Zp[[T]]^* $.
In particular, Conjecture~\ref{mainconjecture}(3) holds for~$ \chi $.
\end{proposition}

\begin{proof}
Conjecture~\ref{muconjecture} holds for $ k = \Q $ by Ferrero and Washington \cite{Fer-Wash}
so we can apply Corollary~\ref{lambdappowercase}. Because $ \l({\bf 1}_\Q) = -1 $
by Lemma~\ref{vlemma}(1), we have $ \l(\chi) = -1 $ as well.
Since $ \chi $ is of type~$ W $, we have
$ \l(H_\chi) = 1 $, hence $ \l(G_\chi) = 0 $. The power
of~$ p $ in $ G_\chi(T) $ is again determined by \cite{Fer-Wash}.
\end{proof}

\begin{corollary} \label{regcor}
If $ p $ is a prime number, then
\cite[Conjecture~3.18]{BBDJR} holds in full for all $ n \ge 2 $
and 1-dimensional Artin characters $ \psi = \psi_\pi : G_\Q \to \Q(\mu_\infty) $
in the following cases:
\begin{itemize}
\item
$ \psi $ is of type $ W $ and $ n \equiv 1 $ modulo the order $ \phi(q) $
of $ \o_p $;

\item
$ \psi = \tilde \psi \eta $ with $ \tilde \psi $
of type~$ W $, $ \eta $ the character of order~2 of $ \Gal(\Q(\mu_q) / \Q) $
viewed as character of $ G_\Q $, and $ n \equiv 1 + \phi(q)/2 $ modulo~$ \phi(q) $.
\end{itemize}
In particular, the $ p $-adic regulator in the conjecture
is a unit in those cases.
\end{corollary}

\begin{proof}
Because of \cite[Proposition~4.17]{BBDJR}, we only need to verify part~(4) of the conjecture.
Note that the characters $ \o_p^{\phi(q)/2} $ and $ \eta^\tau $
coincide for every embedding $ \tau : \Q(\mu_\infty) \to \Qpbar $.
Therefore the $ 1 $-dimensional characters $ \psi $ such that
all $ \psi^\tau \o_p^{1-n} $
for a given~$ n $ are of type $ W $ include the ones given, and one sees easily
using the same train of thought that those are all.
It follows from Lemma~\ref{vlemma} or Proposition~\ref{vprop}
that all resulting $ \Lp p {\Q} {\psi^\tau\o_p^{1-n}} s $ have value a unit at $ n $ (and,
in fact, at every point of its domain where they are defined).
That the regulator in part~(4) of the conjecture is a unit then
follows from part~(2) of the conjecture.
\end{proof}

We determine some $\lambda$-invariants using Corollary~\ref{lambdappowercase}. In all cases where $ k \ne \Q $, Conjecture~\ref{muconjecture} was verified for these characters, see Theorem~\ref{th:proofmu=0}, so the corollary applies. All $ \lambda $-invariants agree with the values found in the computations.

\begin{example} \label{firstexample}
Let $k = \Q$ and $p = 2$. Consider the conductor $3247 = 17 \cdot 191$. Let $\chi$ be a character
of order $16$ and conductor $3247$. We find that $\lambda(\chi)  = 2^2 + 2^4 - 1 = 19$. 
\end{example}

\begin{example} \label{largeexample}
Let $k = \Q$ and $p = 3$. Consider the conductor $2917$. It is a prime number and $2916 = 2^2 \cdot 3^6$. Any character $\chi$ of conductor $2917$ and order $3$, $9$, $27$, $81$, $243$ or $729$ satisfies $\lambda(\chi) = 3^5 - 1 = 242$. This is the largest value of the~$\lambda$-invariant that we found in our computations.
\end{example}

\begin{example}
Let $k = \Q$ and $p = 2$. Consider the conductor $3855 = 3 \cdot 5 \cdot 257$.  Any character $\chi$ of conductor $3855$ and order $4$, $8$, $16$, $32$, $64$, $128$ or $256$ satisfies $\lambda(\chi) = 2^0 + 2^0 + 2^6 - 1 = 65$.
\end{example}

In Examples~\ref{firstquadraticexample} through~\ref{exl1=0}
below, we computed $ \l({\bf 1}_k) $ as in Remark~\ref{compl1}.
With notation as in that remark, we computed in Example~\ref{firstquadraticexample} that $ \l(\chi_k) = -1 $ with $k = \Q(\sqrt{2})$ using
Corollary~\ref{lambdappowercase}, so that $ \l({\bf 1}_k ) = \frac12 (\lambda(\chi_k - 1 ) = -1 $
because $ \Q(\sqrt{2})/\Q $ is the first layer in the 2-cyclotomic
tower of~$ \Q $.

\begin{example} \label{firstquadraticexample}
Let $k = \Q(\sqrt{2})$ and $p = 2$. Consider the conductor $\mathfrak{f} = \mathfrak{p}_{257}$ where~$\mathfrak{p}_{257}$ is one of the two prime ideals above $257$ in $k$. The ray class group of $k$ modulo~$\mathfrak{f}$ has order $2$. Let $\chi$ be the character of conductor $\mathfrak{f}$. We have $k = \Q_1$, the first layer of the cyclotomic $\Z_2$-extension of $\Q$, thus $e = 1$. We find that $\lambda({\bf 1}_k) = -1$, and therefore $\lambda(\chi) = 2^5-1 = 31$.
\end{example}

\begin{example}
Let $k = \Q(\sqrt{5})$ and $p = 5$. Consider the conductor $\mathfrak{f} = \mathfrak{p}_{11}\mathfrak{p}_{31}$ where $\mathfrak{p}_{11}$, respectively $\mathfrak{p}_{31}$,  is one of the two prime ideals above $11$, respectively~$31$, in $k$.  The ray class group of $k$ modulo $\mathfrak{f}$ has order $5$.   Let $\chi$ be a character of conductor $\mathfrak{f}$ and of order $5$. We have~$e = 0$ and $\lambda({\bf 1}_k) = -1$, thus we find that $\lambda(\chi) = 5^0 + 5^0 - 1 = 1$.
 \end{example}

\begin{example}
Let $k = \Q(\sqrt{5})$ and $p = 3$. Consider the conductor $\mathfrak{f} = 2\mathfrak{p}_{109}$ where $\mathfrak{p}_{109}$ is one of the two prime ideals above $109$ in $k$.  The ray class group of~$k$ modulo $\mathfrak{f}$ has order $3$. Let $\chi$ be a character of conductor $\mathfrak{f}$ and order $3$. We have~$e = 0$ and $\lambda({\bf 1}_k) = -1$, and thus $\lambda(\chi) = 3^0 + 3^2 - 1 = 9$.
\end{example}

\begin{example}\label{exl1=0}
Let $k = \Q(\sqrt{33})$ and $p = 2$. Consider the conductor $\mathfrak{f} = \mathfrak{p}_3\mathfrak{p}_{31}$ where $\mathfrak{p}_{3}$ is the only prime ideal above $3$ in $k$ and $\mathfrak{p}_{31}$ is one of the two prime ideals above $31$ in $k$. The ray class group of $k$ modulo $\mathfrak{f}$ has order $2$. Let $\chi$ be the character of conductor $\mathfrak{f}$. We have $e = 0$ and  $\lambda({\bf 1}_k)=0$, thus we find that $\lambda(\chi) = 2^0 + 2^3  = 9$. (This is an example with $\lambda({\bf 1}_k) \geq 0$ as explained in Remark~\ref{compl1}.)
\end{example}

\begin{remark}
The careful reader will have noticed that, in all these examples, the prime numbers, resp.~prime ideals, dividing the conductor but not the prime $p$ are congruent to $1$ modulo $p$, resp.~have their norm congruent to $1$ modulo $p$. This is because they divide exactly the conductor. Indeed, let $\mathfrak{q}$ be a prime dividing~$\mathfrak{f}$ with~$\mathfrak{q}$ not dividing $p$. ($\mathfrak{q}$ and $\mathfrak{f}$ are either rational integers or integral ideals depending on whether $k = \Q$ or a real quadratic field.) Assume that $\mathfrak{f} = \mathfrak{q} \mathfrak{g}$ with $\mathfrak{g}$ not divisible by $\mathfrak{q}$. We have a natural surjection $\mathrm{Cl}_k(\mathfrak{f}) \to \mathrm{Cl}_k(\mathfrak{g})$ whose kernel is a quotient, say~$(\O_k/\mathfrak{q})^\times/H$, of $(\O_k/\mathfrak{q})^\times$. Since $\mathfrak{f}$ is the conductor of $\chi$, this quotient is not in the kernel of the associated Hecke character $\hchi$. But $\hchi$ is non-trivial of~$p$-power order, therefore $(\O_k/\mathfrak{q})^\times/H$ contains an element of order divisible by $p$. It follows that $p$ divides $\mathcal{N}(\mathfrak{q}) - 1$ and thus $\mathcal{N}(\mathfrak{q}) \in 1 + p\Z$. 
\end{remark}

One can also use the theorem of Sinnott when the order of the character is divisible by $p$ but is not a $p$-power. 

\begin{example}
Let $k = \Q$ and $p = 2$. Consider the conductor $889 = 7 \cdot 127$. Let~$\chi$ be a character of conductor $889$ and of order $6$. Write $\chi = \psi\theta$ where $\psi$ has order~$2$ and $\theta$ has order $3$. Then the conductor of $\psi$ is $889$ and the conductor of $\theta$ is $7$. One can check that $\lambda(\theta) = 0$. Since $127 \equiv 1 \mod{7}$, the prime $127$ is totally split in $\Q_\theta$. Therefore $N(\theta, \psi)$ is the number of prime ideals above $127$ in $\Q_\infty$, that is~$2^5 = 32$. In conclusion, we find that $\lambda(\chi) = 32$.
\end{example}

\begin{example}\label{lastexample}
Let $k = \Q$ and $p = 5$. Consider the conductor $1255 = 5 \cdot 251$. Let~$\chi$ be a character of conductor $1255$ and order $10$. (A similar reasoning can be done for the characters of conductor $1255$ and of order $50$ or $250$.) Write $\chi = \theta\psi$ where $\psi$ has order $5$ and $\theta$ has order $2$. As $ \theta $ and $ \psi $ are even, the conductor of $\psi$ is $251$ and the conductor of $\theta=\omega_5^2$ is $5$. One can check that $\lambda(\theta) = 0$. Since $251 \equiv 1 \mod{5}$, the prime $251$ is totally split in $\Q_{\theta\omega_5^3} = \Q(\mu_5)$. Therefore $N(\theta, \psi)$ is the number of prime ideals above $251$ in $\Q_\infty$, that is $5^2 = 25$. In conclusion, we find that $\lambda(\chi) = 25$. 
\end{example}

Now, we go back to the case where $p$ is odd and the order of $\chi$ is prime to $p$. As mentioned in the introduction, our data leads us to make the following conjecture. 

\begin{conjecture} \label{randomconjecture}
Let $p$ be an odd prime and let $d \geq 1$. Define $\mathfrak{X}_d$ to be the set of even $1$-dimensional Artin characters $\chi : G_\Q \to \Qpbar $ of the form~$\chi = \o_p^i\psi$ with $\psi$ of conductor and order prime to $p$, both $\psi$ and $i$ even, and such that~$[\Q_p(\chi):\Q_p] = d$. For $N \geq 1$, let $\mathfrak{X}_d(N)$ be the subset of those characters in $\mathfrak{X}_d$ whose conductor is at
most~$N$. Then we have 
\begin{equation*}
\lim_{N \to +\infty} \frac{\#\{\chi \in \mathfrak{X}_d(N) : \lambda(\chi) > 0\}}{\#\mathfrak{X}_d(N)} = p^{-d}. 
\end{equation*}
\end{conjecture}

To provide some support for the conjecture, we show below graphs for all odd primes up to $37$. Each graph displays two lines, one blue and one red, joining the following points of coordinates $(x,y)$. For both blue and red points, the value of $x$ goes from $0$ to $B$, the upper bound on the conductor (see Table~\ref{tab1}). For the blue points, for a given~$x$, the corresponding value of $y$ is the number of characters in $\mathfrak{X}(x)$ with positive $\lambda$-invariant where $\mathfrak{X}(N) = \bigcup_{d \geq 1} \mathfrak{X}_d(N)$. For the red points, for a given $x$, the corresponding value of $y$ is
\begin{equation*}
\sum_{\chi \in \mathfrak{X}(x)} p^{-[\Q_p(\chi):\Q_p]}.
\end{equation*}
Indeed, under the conjecture, the probability that a random character $\chi $ in $ \mathfrak{X}(x)$ has $\lambda(\chi) > 0 $ tends to $p^{-[\Q_p(\chi):\Q_p]}$ when $x \to \infty$. 

It can be seen on these graphs that the model fits the experimental data quite well for most primes, but is sometimes a bit off for $17$ and $37$, and more so for $19$, $23$ and $29$. This could be due to an incorrect model or simply the fact that our sample set is not large enough since we have some bound on the conductor. Still, the behaviour of both graphs as $ x $ grows seems to generally match well for all primes. 

\begin{figure}[ht]
    \includegraphics[width=\textwidth]{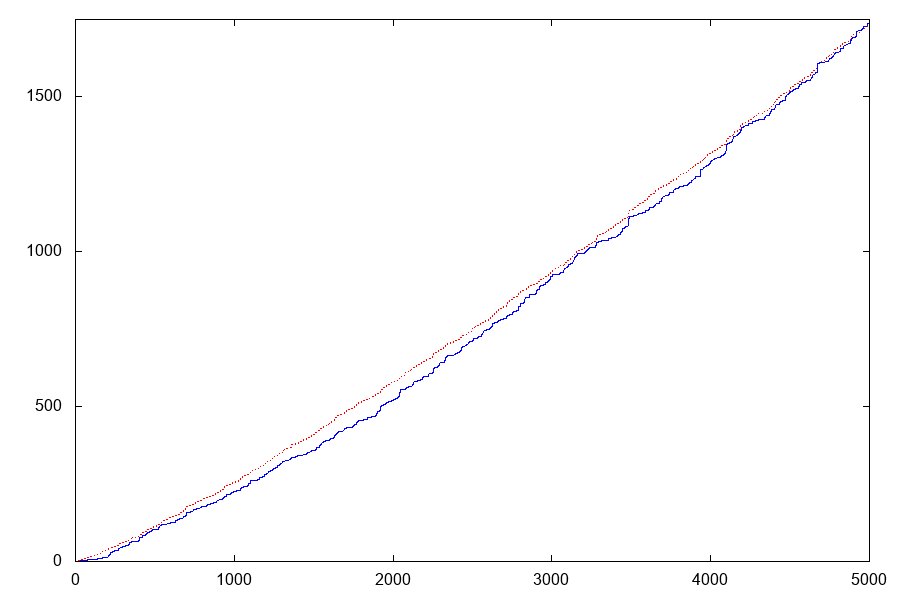}
    \caption*{$p=3$}
\end{figure}

\begin{figure}[ht]
    \includegraphics[width=\textwidth]{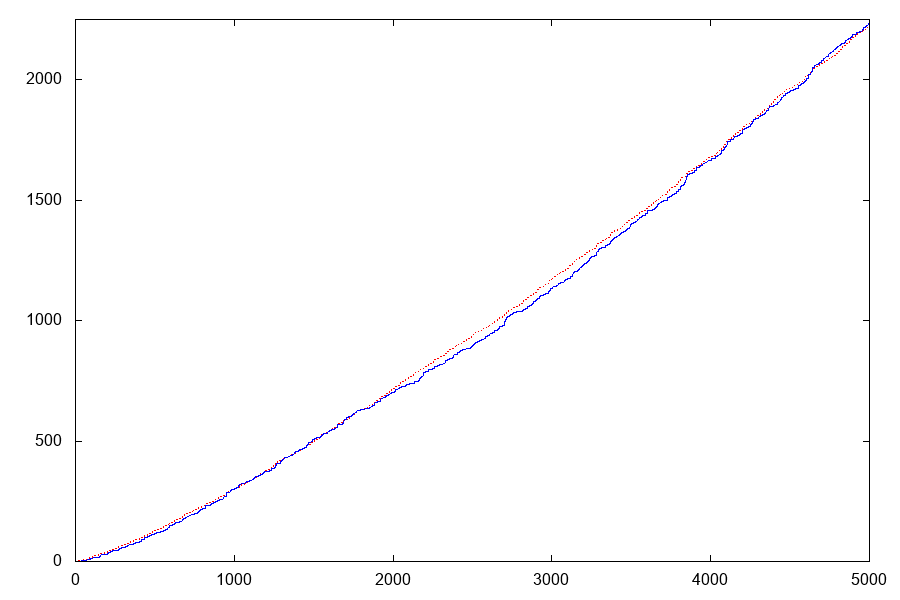}
    \caption*{$p=5$}
\end{figure}

\begin{figure}[ht]
    \includegraphics[width=\textwidth]{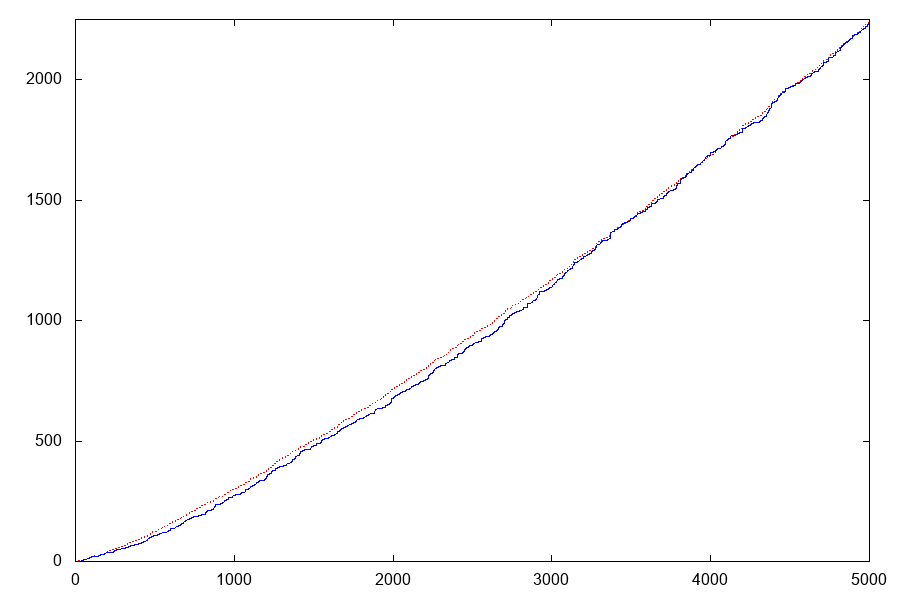}
    \caption*{$p=7$}
\end{figure}

\begin{figure}[ht]
    \includegraphics[width=\textwidth]{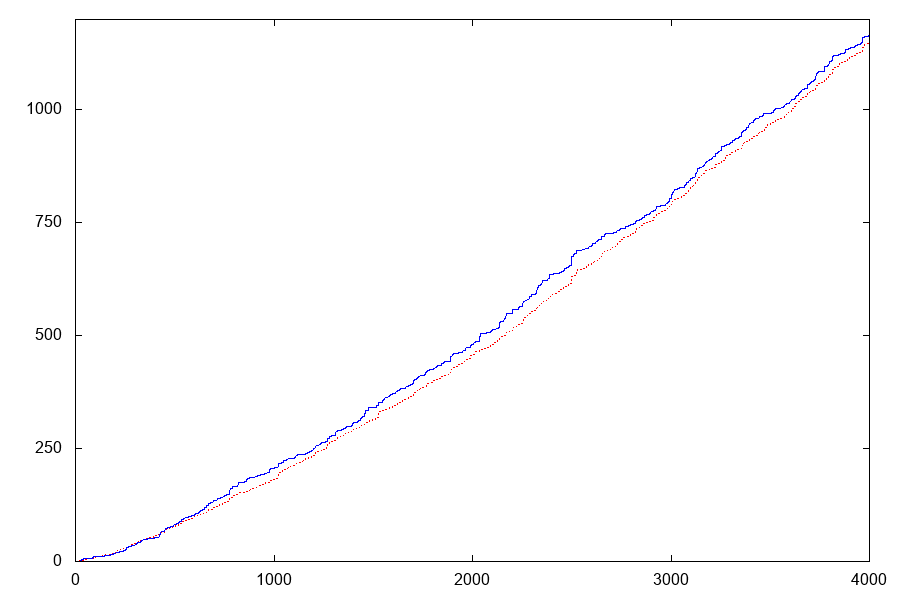}
    \caption*{$p=11$}
\end{figure}

\begin{figure}[ht]
    \includegraphics[width=\textwidth]{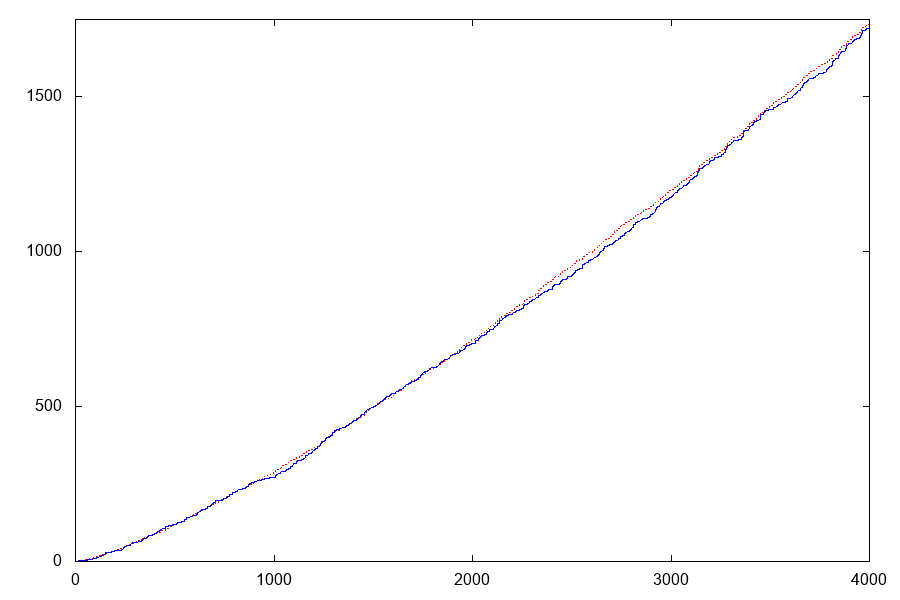}
    \caption*{$p=13$}
\end{figure}

\begin{figure}[ht]
    \includegraphics[width=\textwidth]{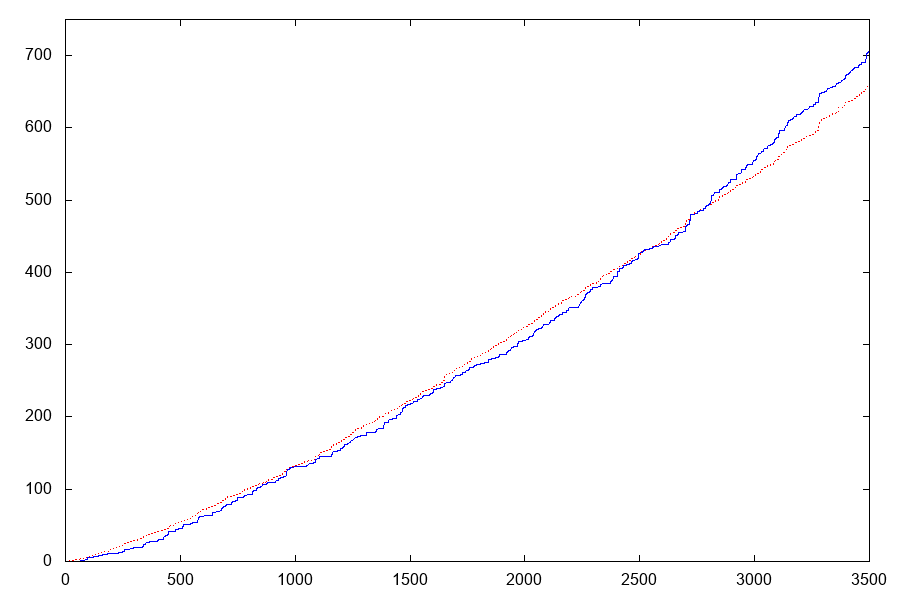}
    \caption*{$p=17$}
\end{figure}

\begin{figure}[ht]
    \includegraphics[width=\textwidth]{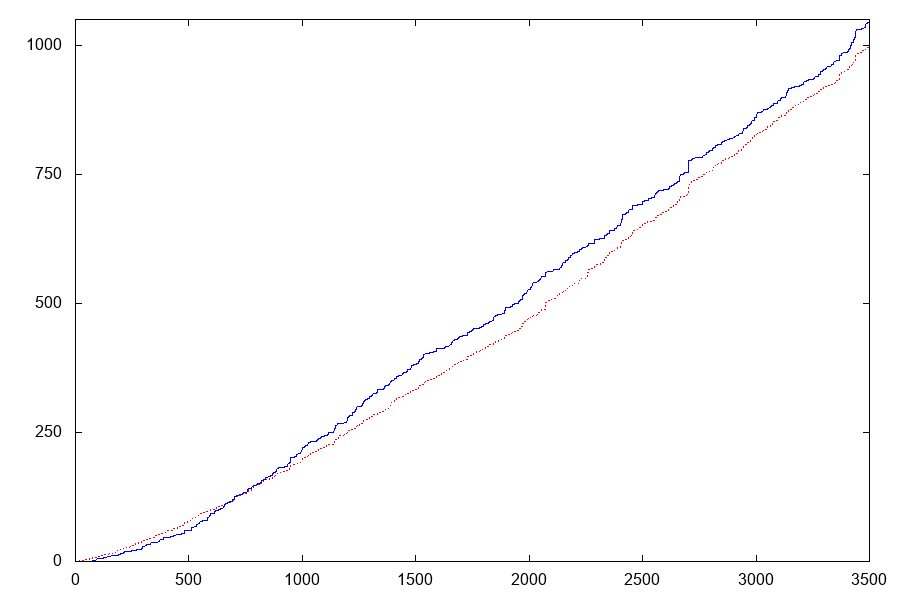}
    \caption*{$p=19$}
\end{figure}

\begin{figure}[ht]
    \includegraphics[width=\textwidth]{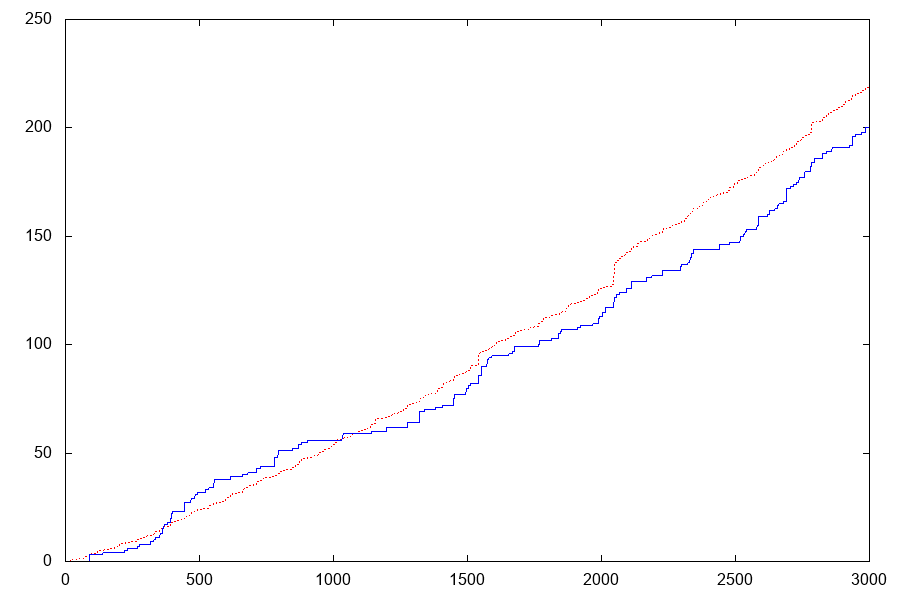}
    \caption*{$p=23$}
\end{figure}

\begin{figure}[ht]
    \includegraphics[width=\textwidth]{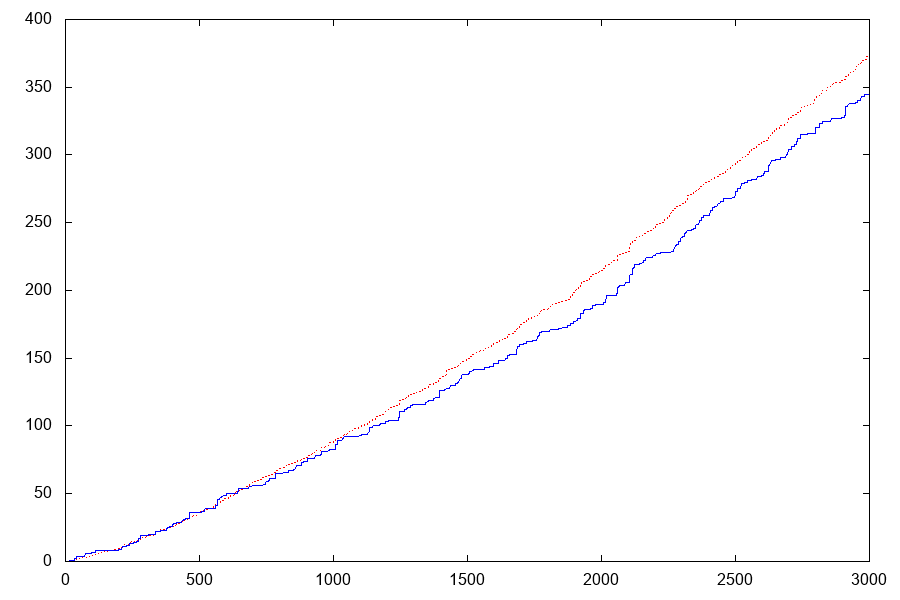}
    \caption*{$p=29$}
\end{figure}

\begin{figure}[ht]
    \includegraphics[width=\textwidth]{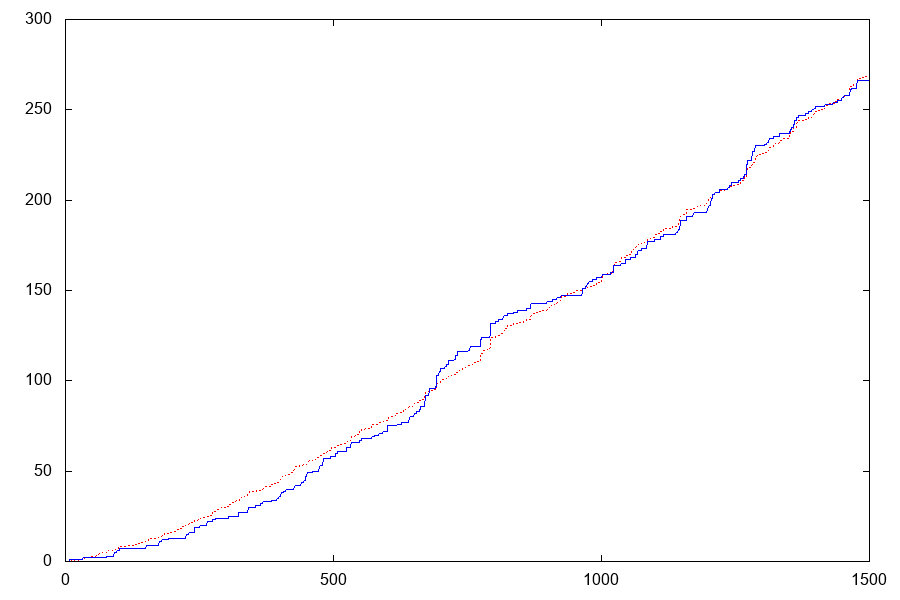}
    \caption*{$p=31$}
\end{figure}

\begin{figure}[ht]
    \includegraphics[width=\textwidth]{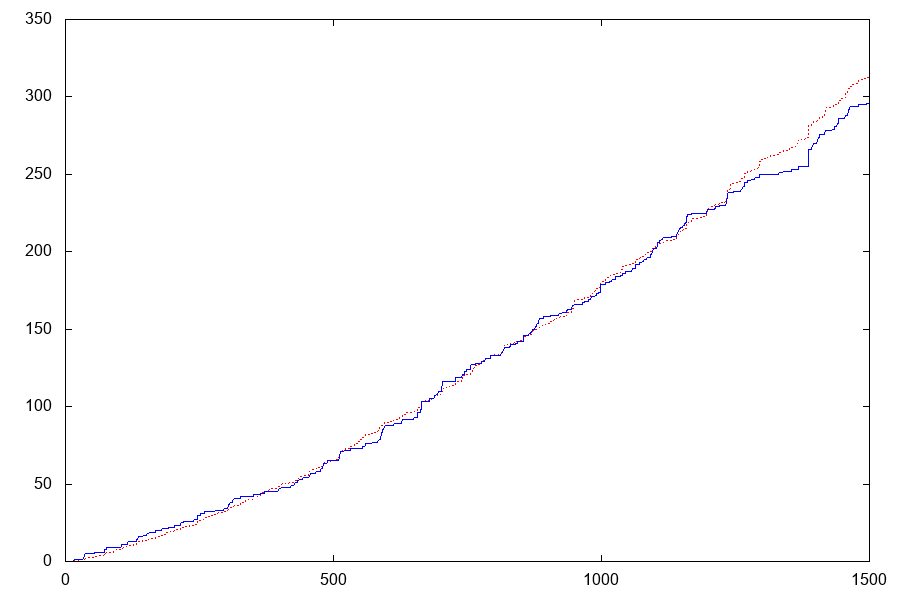}
    \caption*{$p=37$}
\end{figure}